\documentclass[12pt]{elsarticle}

\usepackage{lineno,hyperref}
\modulolinenumbers[5]

\journal{arXiv}
 \usepackage[utf8]{inputenc}
 \usepackage[T1]{fontenc}
 \usepackage{lmodern}
 \usepackage{amsmath}
 \usepackage{amssymb}
 \usepackage{amsthm}
 \usepackage{dsfont}
 \usepackage{graphicx}
 \usepackage{listings}
 \usepackage{subcaption}
 \usepackage{algorithm}
 \usepackage{algpseudocode}
 \usepackage{xspace}
 \usepackage{todonotes}
 \usepackage{csquotes}
 \usepackage{caption}
 \usepackage{nameref}
 \usepackage{cleveref}
 \usepackage{booktabs}
 \usepackage{cancel}
 \usepackage{comment}
 \usepackage{color}
 \usepackage{bm}
 
\theoremstyle{definition}

\theoremstyle{remark}

\def\kg{\text{kg}}
\def\m{\text{m}}
\def\s{\text{s}}

\def\a{\alpha}
\def\b{\beta}
\def\dt{\Delta t}
\def\dx{\Delta x}

\def\Bd{\mathbf{d}}
\def\Bc{\mathbf{c}}
\def\sign{\text{sign}}

\begin{document}

\begin{frontmatter}

\title{Simulation of Gas Mixture Dynamics in a Pipeline Network using Explicit Staggered-Grid Discretization}

%\tnotetext[mytitlenote]{Fully documented templates are available in the elsarticle package on %\href{http://www.ctan.org/tex-archive/macros/latex/contrib/elsarticle}{CTAN}.}

% \author{Yan Brodskyi}
% \address{Humboldt University - Berlin}
% \fntext[footnote1]{Berlin, Germany.}

% \author{Vitaliy Gyrya, Anatoly Zlotnik}
% \address{Los Alamos National Laboratory}
% \fntext[footnote2]{Los Alamos, New Mexico, 87545.}

\author[address1,address2]{Yan Brodskyi}
\author[address2]{Vitaliy Gyrya}
\author[address2]{Anatoly Zlotnik\corref{mycorrespondingauthor}}
\cortext[mycorrespondingauthor]{Corresponding author}
\ead{azlotnik@lanl.gov}


\address[address1]{Humboldt-Universit\"at zu Berlin, Berlin, Germany}
\address[address2]{Los Alamos National Laboratory, Los Alamos, New Mexico, USA}

\begin{abstract}
We develop an explicit second order staggered finite difference discretization scheme for simulating the transport of highly heterogeneous gas mixtures through pipeline networks.  This study is motivated by the proposed blending of hydrogen into natural gas pipelines to reduce end use carbon emissions while using existing pipeline systems throughout their planned lifetimes.  Our computational method accommodates an arbitrary number of constituent gases with very different physical properties that may be injected into a network with significant spatiotemporal variation.  In this setting, the gas flow physics are highly location- and time- dependent, so that local composition and nodal mixing must be accounted for.  The resulting conservation laws are formulated in terms of pressure, partial densities and flows, and volumetric and mass fractions of the constituents.  We include non-ideal equations of state that employ linear approximations of gas compressibility factors, so that the pressure dynamics propagate locally according to a variable wave speed that depends on mixture composition and density.  We derive compatibility relationships for network edge domain boundary values that are significantly more complex than in the case of a homogeneous gas.  The simulation method is evaluated on initial boundary value problems for a single pipe and a small network, is cross-validated with a lumped element simulation, and used to demonstrate a local monitoring and control policy for maintaining allowable concentration levels.
\end{abstract}


\end{frontmatter}

%\newpage 

%\tableofcontents

%\newpage

%=======================================================

\section{Introduction} \label{sec:intro}

%=======================================================

% 1 - background

The accelerating transition of industrial economies from reliance on fossil fuels to renewable sources of energy has motivated the development of new methods for production, transport, and utilization of hydrogen \cite{yue2021hydrogen}.   Hydrogen does not create carbon dioxide or other pollutants when burned or used in a fuel cell, and can be produced by electrolysis using sustainable electricity sources \cite{dawood2020hydrogen}.  It has been examined as an additive to natural gas for general end use \cite{nitschke2012admixture}, and can be delivered through existing gas pipeline infrastructure.  The blending of hydrogen into natural gas pipelines creates engineering and operating issues \cite{melaina2013blending}, and new mathematical approaches are being developed to extend pipeline flow modeling and simulation methods to heterogeneous gas mixtures \cite{uilhoorn2009dynamic}.

% 2 - Motivation

The introduction of hydrogen into natural gas transport systems raises significant issues related to materials and engineering \cite{melaina2013blending,erdener2023review}. Hydrogen can cause embrittlement of certain types of steel, which can affect pipeline integrity and the reliability of compressor turbomachinery \cite{hafsi2019computational}. Supposing that such technological issues are addressed, the significant differences in physical and chemical properties of hydrogen from those of natural gas will cause the pipeline flow dynamics and energy released in combustion of a blend of these gases to vary broadly depending on its composition \cite{abd2021evaluation}.   The assumption of homogeneous composition is generally sufficient to predict the behavior of a pipeline system that transports processed natural gas to consumers with time-varying demands.  For a natural gas pipeline that accepts time-varying hydrogen injections at multiple locations, and whose end users consume gas with time-varying profiles, it is critical for the variation in composition to be incorporated in predictive simulation as well as real-time monitoring \cite{erickson2022importance}.  This will ensure appropriate modeling of physical flow, enable accounting for the form in which energy is actually delivered, quantify the effects of blending on pipeline efficiency, and facilitate assessment of changes to downstream carbon emissions \cite{zlotnik2023effects}.

% 3 - Issues and requirements

The primary motivation for pipeline simulation that accounts for blending of heterogeneous gases is to evaluate the effects of such modifications on energy transmission capacity.  Because the properties of hydrogen are so different from those of natural gas, modeling the flow physics and energy content of blends of these gases leads to mathematical settings that are multi-scale, numerically challenging, and highly sensitive to mass fractions in both the steady-state and transient flow regimes.  Accurate partial differential equation (PDE) representations of pressures, flows, and calorific values are critical when simulating blending of heterogeneous gases through complex, large-scale pipeline networks.  Modeling the flow of a homogeneous gas on a network requires PDEs for mass and momentum conservation on each pipe, as well as a linear mass flow balance equation and one boundary condition at each network junction. For each additional constituent gas, the modeling requires another PDE on each pipe, another nonlinear nodal balance equation at each junction, and another boundary condition parameter at each node where gases enter the network.  Additional corresponding state variables are needed to account for changes in mass fraction, which affect total density,  flow dynamics, and energy content.  In the case of hydrogen, the faster wave speed corresponding to lower density aggravates the numerical ill-conditioning of the associated dynamic model.

% 4 - Previous studies

New modeling approaches have been under development over the past decades in order to characterize the phenomena that result from hydrogen blending starting with a study of modifications needed for simulation \cite{uilhoorn2009dynamic}. An early study examined the problem of fuel minimization of  hydrogen-methane blends and noted the trade-offs between delivery pressure, hydrogen fraction, and transmitted energy for a simple pipeline \cite{tabkhi2008mathematical}.  Others have examined various aspects of this simulation problem \cite{guandalini2017dynamic, fan2021transient}, including with the use of reduced-order modeling \cite{agaie2017reduced}.  Conditions under which pipeline pressures may exceed allowable upper limits were explored, and it was shown that the likelihood of this occurrence increases proportionally with increasing hydrogen concentration \cite{hafsi2019computational}.  Another study examined the effects of hydrogen blending on methods for detecting, localizing, and estimating leaks \cite{subani2017leak}, and showed that the leak discharge is expected to increase with hydrogen concentration.  A moving grid method and an implicit backward difference method for tracking gas concentration were both shown to perform well, but the implicit difference method was observed to lose accuracy because of numerical diffusion \cite{chaczykowski2018gas}.   A numerical simulation scheme for transient flows on cyclic networks with homogeneous blends was developed using the method of characteristics \cite{elaoud2017numerical}.  Another study modeled composition tracking in pipeline networks \cite{witkowski2017comprehensive}, without incorporating the control actions of compressors.  

The  models in the studies described above demonstrate a simulation capability or conduct a sensitivity study for a specific network. A recent study has examined the well-known challenges of optimizing flows of mixing gases throughout a network, which leads to highly challenging nonlinear mixed-integer programming formulations \cite{hante2019complementarity,hante2023gas}.  Addressing the challenging conceptual questions related to design, operational, and economic issues of pipeline transport of gas mixtures requires minimal and generalizable mathematical models that adequately describe the flow physics, in addition to more complex frameworks that comprehensively characterize all flow properties.  Studies that address generalizable network modeling for use in simulation and optimal control of hydrogen blending in natural gas pipelines have only emerged in recent years \cite{zhang2022modelling}.

% 5 - Challenges and gaps

Several conceptual and computational challenges remain open in order to characterize the effect of hydrogen blending on gas pipeline transients.  The more physically complex flows on each pipe in the network must be accurately resolved, and nodal compatibility conditions are needed to represent the distributed dynamic flow of mixtures of gases through a network with time-varying injections and withdrawals of heterogeneous constituents, as well as compressor controls.  The set memberships of pipes with physical flows that are incoming and outgoing at a node are needed to appropriately evaluate the mixing conditions.  These set memberships change with flow reversals, and this leads to mixed-integer programs that are challenging in steady-state \cite{kazi2024modeling}. In transient simulation, it is standard to assume that flows do not reverse during the simulation period \cite{chaczykowski2018gas}. 

% 6 - This study

In this study, we develop an explicit second order staggered finite difference discretization scheme for solving initial boundary value problems (IBVPs) for simulating the transport of highly heterogeneous gas mixtures through pipeline networks.  The approach is based on a recent method for simulating the flow of a homogeneous gas throughout a pipeline network \cite{gyrya2019explicit}, and inherits its desirable properties of stability and second-order accuracy.  The  method  described here accommodates an arbitrary number of constituent gases that may have very different physical properties.  The boundary conditions can include injection of pure or mixed gases into the network with significant spatio-temporal variation and time-varying withdrawal.  The conservation laws are formulated in terms of pressure, partial densities and flows, and volumetric and mass fractions of the constituents.  Our modeling includes non-ideal equations of state that we develop using linear approximation of gas compressibility factors, so that the pressure dynamics propagate locally according to a variable wave speed that depends on mixture composition and density.  We also derive compatibility relationships for network nodes, which may include compressors, that are significantly more complex than in the case of a homogeneous gas.  The simulation method is evaluated on initial boundary value problems for a single pipe and a small network, is  cross-validated with a lumped element simulation, and used to demonstrate a local monitoring and control policy for maintaining allowable concentration levels.

% 7 - Manuscript structure

The rest of the manuscript is structured as follows.  In Section \ref{sec:modeling}, we describe the physical modeling of gas mixture transport for a single pipe, our approach to non-ideal equation of state modeling for a gas with variable composition of multiple constituents with different properties, and the treatment of boundary conditions on a pipeline network.  Section \ref{sec:discretization} contains our major contribution, in which we describe the staggered-grid discretization for the system of PDEs on a network, which is used to compute the evolution of the pressure, flow, and concentration variables in each pipe and at nodes in the network.  Then, Section \ref{sec:monitoring} describes a monitoring and corrective control policy that can be used to ensure that mass fractions of constituent gases are maintained within acceptable limits.  Section \ref{sec:computational} contains a collection of computational studies to address key physical modeling and numerical analysis issues, to compare modeling of ideal and non-ideal gas flows, to demonstrate the simulation method for a small test network system, and cross-validate the method with another published approach for gas mixtures or a homogeneous gas.  Finally, in Section \ref{sec:conclusions} we conclude with a discussion of the outcomes of this study, connections with contemporary work, and future directions.

%=======================================================

\section{Modeling of Gas Pipeline Flows} \label{sec:modeling}

%=======================================================

There is a rich body of literature on modeling and simulation of natural gas pipeline flows \cite{wylie1978fluid,thorley1987unsteady,osiadacz1984simulation}.  Our goal in this study is an explicit numerical method for simulating flows on a large-scale pipeline network with heterogeneous gas injections.  Based on the approach taken in our previous study \cite{gyrya2019explicit}, we model a pipeline network as a set of edges to represent pipes that are connected at nodes that represent junctions.  The network topology is a connected directed metric graph $(\mathcal{V},\mathcal{E})$ where $\mathcal{V}$ and $\mathcal{E}$ are used to denote sets of nodes and edges, respectively.  An edge $(i,j)\subset \mathcal{E}$ connects two nodes $i,\,j\in\mathcal{V}$.  The dynamics of gas flows on the whole network are determined by the gas flow physics on each edge, as well as compatibility and boundary conditions that are defined on nodes. In this section, we first describe physical flow modeling on a single pipe, including equation of state computation used to determine pressures based on the local composition of an arbitrary number of mixing gases.  We then describe the notation and treatment of boundary conditions over the network.

%=======================================================

\subsection{Gas Mixture Flow in a Pipe} \label{sec:onepipe}

%=======================================================

The established model of isothermal homogeneous gas flowing through a single pipe is described using the Euler equations in one dimension,
\begin{subequations} \label{eq:gaspde0}
\begin{align}
    \frac{\partial d}{\partial t} + \frac{\partial (d v)}{\partial x} & = 0, \label{eq:gaspde0a} \\
    \frac{\partial (d v)}{\partial t} + \frac{\partial (d v^2 + p)}{\partial x} & = - \frac{\lambda}{2D}d v |v| - d g \frac{\partial h}{\partial x}, \label{eq:gaspde0b} \\
    p = Z(p,T)RT \cdot d & = \sigma^2 \cdot d. \label{eq:gaspde0c}
\end{align}
\end{subequations}
Equations \eqref{eq:gaspde0} represent mass conservation, momentum conservation, and the gas equation of state law.  Here the state variables $u$, $p$, and $d$ represent gas velocity, pressure, and density, respectively, with dependence on time $t$ and space $x\in[0,L]$, where $L$ is the length of the pipe, and the variable $h$ gives the elevation of the pipeline.  The dimensionless parameter $\lambda$ is a factor in the phenomenological Darcy-Weisbach term that models momentum loss caused by turbulent friction.  The other parameters are the internal pipe diameter $D$, and the wave (sound) speed  $\sigma=\sqrt{Z(p,T)RT}$ in the gas where $Z(p,T)$, $R$, and $T$ are the gas compressibility factor, specific gas constant, and absolute temperature, respectively, and the gravitational acceleration constant $g$.  We define the per-area mass flux as $\phi=vd$, and apply baseline assumptions for gas transmission pipelines for simplicity and to focus on the aspects of mixing dynamics.  We suppose that each pipe is horizontal and has uniform diameter and internal surface roughness, that the flow is turbulent and has high Reynolds number, that gas flow is an isothermal process at constant and uniform temperature and is adiabatic so that there is no heat exchange with the ground.  The gas compressibility factor $Z$ depends on pressure and temperature, and in what follows will also depend on gas composition.  

Suppose now that the flowing gas consists of $n$ gas components, and let $d^\a$ denote the partial densities (measured in $\kg$/$\m^3$) of the gas components with $\a=1,\dots,n$, which sum to the total density according to $d = d_1+\dots+d_n$.  In addition, let $c^\a=d^\a/d$ denote the mass fractions of the gas components, which sum to unity according to $1 = c^1+\dots+c^n$.  We use the notation $\Bd:=(d^1,\dots,d^n)$ and $\Bc:=(c^1,\dots,c^n)$ to denote the vectors of partial densities and mass fractions of the gas components.  We may henceforth use the term concentration to refer to the mass fraction, and where the volumetric fraction is discussed then this will be referred to explicitly. Under the above assumptions, the conservation equations \eqref{eq:gaspde0a}-\eqref{eq:gaspde0b} can be extended to conservation laws for flow of the gas mixture in the form
\begin{subequations} \label{eq:multigaspde0}
\begin{align}
    \frac{\partial}{\partial t} d^\a + \frac{\partial}{\partial x}(d^\a v) & = \epsilon^\a \Delta d^\a, \quad \text{for all}\  \a=1,\dots, n, \label{eq:multigaspde0a} \\
    \frac{\partial}{\partial t}(d v) + \frac{\partial }{\partial x} (d v^2 + p) & = - \frac{\lambda}{2D}d v |v|. \label{eq:multigaspde0b} 
\end{align}
\end{subequations}
The explicit dependence of the pressure on the gas component mass fractions $c^1,\ldots,c^n$ will be discussed in detail in Section \ref{sec:EOS}.  The parameters $\epsilon^\a$ (cm$^2$/s) are diffusion coefficients for individual gas components, and $\Delta$ denotes the Laplacian operator in equation \eqref{eq:multigaspde0a}.  We proceed to omit the diffusion terms on the right hand side of the conservation of mass equation \eqref{eq:multigaspde0a} because relatively small diffusion coefficients $\epsilon^\a$ result in a negligible influence on advection dynamics, which will be confirmed later in the computational results in Section \ref{sec:onepipe_example}.  We also omit the term $\partial (d v^2)/\partial x$ that quantifies the effects of kinetic energy on momentum conservation in \eqref{eq:multigaspde0b}, as it is several orders of magnitude lower than the other terms in equation \eqref{eq:gaspde0b} in the regime of interest \cite{osiadacz1984simulation}. Using the above simplifications, the equations \eqref{eq:multigaspde0a}-\eqref{eq:multigaspde0b} can be written in terms of (per area) mass flux $\phi = d v$ (in units of $\kg/\m^{2}/\s$) as 
\begin{subequations} \label{eq:multigaspde1}
\begin{align}
    \frac{\partial}{\partial t} d^\a + \frac{\partial}{\partial x}(c^\a \phi) & = 0, \quad \text{for all}\  \a=1,\dots, n, \label{eq:multigaspde1a} \\
    \frac{\partial}{\partial t}\phi + \frac{\partial }{\partial x} p & = - \frac{\lambda}{2D}\frac{|\phi|\phi}{d}. \label{eq:multigaspde1b} 
\end{align}
\end{subequations}
We now examine how to modify the equation of state \eqref{eq:gaspde0c} to facilitate the treatment of gas mixtures.

%=======================================================

\subsection{Equation of State Modeling} \label{sec:EOS}

%=======================================================

The modeling and implementation of equations of state in gas pipeline simulation is a complex topic that has seen significant development for many decades \cite{lee1966viscosity,peng1976new,mccarty1974modified}, which has culminated in the widely used AGA8 \cite{starling1992compressibility} and GERG \cite{kunz2012gerg} calculators.  The use of the latter two approaches has been compared in a recent study on the potential impacts of hydrogen blending \cite{bainier2019impacts}.  Such implementations require an inherently implicit approach to computation, which we wish to avoid in this study where we aim for a fully explicit simulation method. 

To establish a setting for our treatment of the state equation, we define additional notation here.  We have already defined the \emph{total density} $d$, the \emph{partial density} $d^\a$ of each gas component, the \emph{mass fraction} $c^\a=d^\a/d$ of each gas component, and the \emph{total mass flux} $\phi= v d$ of the mixture in Section \ref{sec:onepipe}.  We further define the \emph{partial mass flux} $\phi^\a = v d^\a$, where $v$ is the \emph{flow velocity}, and define the \emph{volumetric fraction} $\gamma^\a$ of the gas as the fraction of a unit volume occupied by gas component $\alpha$. The volumetric fractions of the gas mixture satisfy $1=\gamma^1 + \cdots + \gamma^n$.  Finally, we define the \emph{individual density} $\rho^\a$ of each gas component, which refers to the density of gas component $\a$ within the fraction of the volume that it occupies in mixture.  In contrast to the partial densities $d^\a$, the individual densities $\rho^a$ are not additive (e.g., $\rho^1=\rho^2=2d^\a=d^\a/\gamma^\a=d$ in an equal mixture of two gases $\a=1,2$ with the same properties).  

To expand on the distinction between partial densities and individual densities, suppose that each constituent gas $\alpha\in\{1,\ldots,n\}$ is known to satisfy a \emph{pressure-density mapping} $\pi_\a$, which relates pressure $p^\a$ and individual density $\rho^\a$ for the gas $\alpha$ according to
\begin{align}\label{eq:pressure-density_maps}
	p^\a=\pi_\alpha(\rho^\a)  \quad \text{and} \quad  \rho^\a=\pi_\alpha^{-1}(p^\a),
\end{align}
where $\pi_\a$ and $\pi_\a^{-1}$ are known, possibly tabulated functions.  We assume that $\pi_\a$ is a continuous and monotonically strictly increasing function for each gas $\alpha$ for constant temperature $T$, so that $\pi_\a$ and thus $\pi_\a^{-1}$ are bijective.  Using the pressure-density mapping concept, we see that the pressure for all gas components in a mixture should be the same, i.e.,
\begin{align}\label{eq:pressure_equivalence}
	\pi_\a(\rho^\a)=p^\a=p=p^\b=\pi_\b(\rho^\b) \qquad \text{for all}\  \, \a,\b\in\{1,\dots,n\}.
\end{align}
Using the notations defined above, we can relate the partial density, individual density, and volumetric fraction of each gas component $\a$ according to
\begin{align}\label{eq:density_relations}
	d^\a = \gamma^\a \rho^\a ,\qquad 	\rho^\a = d^\a/\gamma^\a,	\qquad \text{and} \qquad 	\gamma^\a = d^\a/\rho^\a.
\end{align}
We re-iterate also that the partial densities, volumetric fractions, and mass fractions of the gas components satisfy
\begin{align}\label{eq:fraction_sums}
	d = \sum_{\a=1}^n d^\a, 	\qquad \sum_{\a=1}^n \gamma^\a = 1, \qquad \text{and} \qquad \sum_{\a=1}^n c^\a = 1.
\end{align}
The pressure balance \eqref{eq:pressure_equivalence} together with density conversion relations \eqref{eq:density_relations} and 
volume fraction condition \eqref{eq:fraction_sums} yields the following condition on pressure:
\begin{align}\label{eq:EOS_for_mixture}
	\sum_\a \gamma^\a 	= 	1	\quad \text{and} \quad 	\gamma^\a = \frac{d^\a}{\rho^\a} 	\qquad \implies \qquad
	\sum_\a \frac{d^\a}{\pi_\a^{-1}(p)} = \sum_\a \frac{d^\a}{\rho^\a} = 1,
\end{align}
where $\rho^\a=\pi_\a^{-1}(p)$ are the gas-specific density-pressure functions that have been explained in equations
\eqref{eq:pressure-density_maps} and \eqref{eq:pressure_equivalence}.  Equation \eqref{eq:EOS_for_mixture} provides an implicit equation for the pressure $p$ in terms of the partial densities $d^\a$ of the mixture.

Existence and uniqueness of the solution of \eqref{eq:EOS_for_mixture} for pressure follows from the assumptions of continuity and bijectivity of each of the pressure-density mappings $\pi_\a$ in equation \eqref{eq:pressure-density_maps}.  Indeed, assume that $d^\a$ are fixed. In the limit of large pressure $p$, each of the densities $\rho^\a=\pi_\a^{-1}(p)$ goes to infinity and, therefore, $\sum_\a d^\a/\pi_\a^{-1}(p)$ goes to zero.  On the other hand, in the limit of small pressure $p$, each of the densities $\rho_\a=\pi_\a^{-1}(p)$ goes to zero and therefore $\sum_\a d^\a/\pi_\a^{-1}(p)$ goes to infinity.  By the continuity and bijectivity of the functions $\pi_\a$, it follows that there exists a value of pressure $p\in(0,\infty)$ for which \eqref{eq:EOS_for_mixture} holds.

Using the above definitions, we may suppose that the mapping $\pi_\a^{-1}$ can be defined as 
\begin{align} \label{eq:individual_EOS}
	\rho_\alpha=\pi_\a^{-1}(p) = \frac{p}{R_\a T Z_\a(p,T)},
\end{align}
where $Z_\a(p,T)$ is a compressibility factor at the local pressure $p$ and temperature $T$ for the gas $\a$, and $R_\a$ is the specific gas constant.  Further, consider a linear approximation for gas compressibility of form 
\begin{align} \label{eq:EOS_linear_Z}
	Z_\a(p,T) = 1 + a_\a(T) p.
\end{align}
Note that an ideal gas approximation corresponds to taking $a_\a\equiv0$. We will subsequently use linear approximations for the compressibility factors of natural gas and hydrogen in our application-oriented computational studies, with parameters $a_{NG}$ and $a_{H}$, respectively. Explicit calculations of these coefficients based on experimental tables \cite{kareem2016compng,mihara1977comph2} can be found in the Appendix \ref{appendix:compressibility}.

Substituting \eqref{eq:individual_EOS} into \eqref{eq:EOS_for_mixture} and using the linear approximation  \eqref{eq:EOS_linear_Z} for compressibility yields
\begin{align} \label{eq:EOS_linear_compressibility}
	1 	=  \sum_\a \frac{d^\a}{\pi_\a^{-1}(p)} 	= \sum_\a \frac{d^\a R_\a T Z_\a(p)}{p}
	= 	\sum_\a \frac{d^\a R_\a T (1+a_\a(T) p)}{p}.
\end{align}
Multiplying \eqref{eq:EOS_linear_compressibility} by $p$ and solving the resulting linear equation for $p$ yields
an explicit expression for pressure in terms of the partial densities $\Bd=(d^1,\ldots,d^n)$ and the temperature $T$ of form
\begin{align} \label{eq:explicit_pressure}
	p(\Bd) = \frac{\sum_\a d^\a R_\a T}{1 - \sum_\a d^\a R_\a T a_\a(T)}.
\end{align}
Throughout the rest of our study, we suppose for simplicity that temperature is constant and uniform throughout the pipeline network through the duration of a simulation, and thus we can omit the dependence on temperature in equations \eqref{eq:EOS_linear_Z} and \eqref{eq:explicit_pressure}.  The equation of state can then be accounted for in the conservation laws \eqref{eq:multigaspde1} by using equation \eqref{eq:explicit_pressure} in equation \eqref{eq:multigaspde1b} and using additivity of partial densities in equation \eqref{eq:multigaspde1a}, yielding the system of equations
\begin{subequations} \label{eq:multigaspde2}
\begin{align}
    \frac{\partial}{\partial t} d^\a + \frac{\partial}{\partial x}\left(\frac{d^\a}{\sum_\a d^\a} \phi\right) & = 0, \quad \text{for all}\  \a=1,\dots, n, \label{eq:multigaspde2a} \\
    \frac{\partial}{\partial t}\phi + \frac{\partial }{\partial x} p(\Bd) & = - \frac{\lambda}{2D}\frac{|\phi|\phi}{\sum_\a d^\a}, \label{eq:multigaspde2b}  \\
    p(\Bd) & = \frac{\sum_\a d^\a R_\a T}{1 - \sum_\a d^\a R_\a T a_\a(T)},\label{eq:multigaspde2c}
\end{align}
\end{subequations}
which represents the flow dynamics of the mixture on one pipe in terms of total mass flux $\phi$ and partial densities $d^\a$.

%=======================================================

\subsection{Network Modeling and Boundary Conditions} \label{sec:network}

%=======================================================

As noted above, a network of pipelines can be modeled as a directed graph $\mathcal{G}=(\mathcal{V},\mathcal{E})$ where each edge represents a pipe, and the edge metric gives the pipe's physical properties. The set of directed edges $\mathcal{E} \subset \mathcal{V} \times \mathcal{V}$ has elements $(i,j) = k \in \mathcal{E}$, which denotes that pipe $k$ connects the nodes $i,j \in \mathcal{V}$.  Every edge $k \in \mathcal{E}$ is associated with a diameter $D_{k}$, length $L_{k}$, cross-section area $S_{k}$, and friction coefficient $\lambda_{k}$. Gas flow through each pipe $k \in \mathcal{E}$ is characterized by the hydrodynamic equations \eqref{eq:multigaspde2}, where $\phi_k$ and $d_k^\a$ denote flows and partial densities that depend on time $t\geq 0$ and distance $x=[0,L_k]$.  For each node $q  \in \mathcal{V}$, we define two sets of incoming and outgoing pipes by $\partial^{+}{q} = \{ k \in \mathcal{E} ~|~ \exists i \in \mathcal{V} \mbox{ s.t. } k = (i,q)\}$ and by $\partial^{-}{q} = \{ k \in \mathcal{E} ~|~ \exists i \in \mathcal{V} \mbox{ s.t. } k = (q,i)\}$, respectively.   We will use the index $k$ to enumerate the pipes adjacent to a node and append it as a subscript to indicate that a particular quantity corresponds to the pipe $k$.  For each pipe $k\in\mathcal{E}$, the indexed flow equations and equation of state are then given by
\begin{subequations} \label{eq:multigaspde3}
\begin{align}
    \frac{\partial}{\partial t} d_k^\a + \frac{\partial}{\partial x}\left(\frac{d_k^\a}{\sum_\a d_k^\a} \phi_k\right) & = 0, \qquad \text{for all} \ \ \a=1,\dots, n, \label{eq:multigaspde3a} \\
    \frac{\partial}{\partial t}\phi_k + \frac{\partial }{\partial x} p(\Bd_k) & = - \frac{\lambda}{2D}\frac{|\phi_k|\phi_k}{\sum_\a d_k^\a}, \label{eq:multigaspde3b}  \\
    p(\Bd_k) & = \frac{\sum_\a d_k^\a R_\a T}{1 - \sum_\a d_k^\a R_\a T a_\a(T)}, \label{eq:multigaspde3c}
\end{align}
\end{subequations}
where $\Bd_k=(d_k^1,\ldots,d_k^n)$, and where $\phi_k$ and $d_k^\a$ for all $\alpha=1,\ldots,n$ are defined on $t\geq 0$ and for $x\in[0,L_k]$.
We use a shorthand to denote the boundary values of pressure, densities and concentrations of gas $\alpha$, or mass flux on a pipe according to $p_k (t, 0) = \underline{p}_k(t)$ and $p_k (t, L_k) = \overline{p}_k(t)$, $d_k^{\a} (t, 0) = \underline{d}_k^{\a}(t)$ and $d_k^{\alpha} (t, L_k) = \overline{d}_k^{\alpha}(t)$, $c_k^{\alpha} (t, 0) = \underline{c}_k^{\alpha}(t)$ and $c_k^{\alpha} (t, L_k) = \overline{c}_k^{\alpha}(t)$, and $\phi_k (t, 0) = \underline{\phi}_k(t)$ and $\phi_k (t, L_k) = \overline{\phi}_k(t)$.  For auxiliary nodal values of pressure, density, withdrawal flow, injection flow, injection supply concentration, and outflow concentrations at a node $q\in\mathcal V$, we use the notation $p_q(t)$, $\rho_q(t)$, $F_q^d(t)$, $F_q^s(t)$, $c_q^{\alpha,s}$, and $c_q^{\alpha}$, where the latter denotes the concentration for outgoing flows after mixing at the node. In this study, we suppose for the purpose of clearer exposition that the direction of physical gas flow is in the positive oriented direction for each pipe, and that boundary conditions for the considered initial boundary value problems do not result in changes in flow direction.  The developed method can account for changing flow directions with minor modifications to nodal balance conditions.  We may omit the dependence on time and space in the subsequent exposition where this dependence is understood, in order to simplify notation.

The boundary conditions at a node, which represents a junction that connects two or more pipes, are formulated in terms of 
\textit{(i)} conservation of mass for each of the gas components, i.e., flow balance conditions; and
\textit{(ii)} a form of a continuity condition.  The conservation of mass is stated in terms of the flow balance condition at the node.  With respect to a set of pipes oriented as entering and leaving node $q$, the flow balance condition takes the form
\begin{equation}\label{eq:flow_balance_general_form}
    \sum_{k\in\partial^+q} S_k \overline{c}_k^{\a} \overline{\phi}_k 
    - 
    \sum_{k\in\partial^-q} S_k \underline{c}_k^{\a} \underline{\phi}_k =  c_q^{\alpha} F_q^d - c_q^{\alpha,s} F_q^s, \quad \forall \,\alpha\in\{1,\ldots,n\}, \, \forall \, q\in\mathcal{V},
\end{equation}
where we use $c_q^{\alpha,s}$ to denote the concentration of gas $\a$ in an injection flow $F_q^s$, and $c_q^{\alpha}$ denotes the nodal concentration after mixing that is also the concentration of the withdrawal flow $F_q^d$ from node $q$, and $S_k = \pi(D_k/2)^2$ is the cross-sectional area of the pipe $k$ with diameter $D_k$.  We suppose that only one of $F_q^s$ or $F_q^d$ is nonzero at any given time, and the total withdrawal flow from node $q$ can be stated as $F_q=F_q^d-F_q^s$.  Note that the nodal flows $F_q^s$ and $F_q^d$ are total mass flows, whereas the pipe flows $\overline{\phi}_k$ and $\underline{\phi}_k$ give per area mass flux at the pipe endpoints.

The continuity condition takes a non-trivial form for two reasons.  First, we are neglecting the diffusive term in the conservation of momentum equation, for which we provide computational justification in Section \ref{sec:onepipe_example}.  
Second, we model junctions as points with complete and instantaneous mixing of incoming and outgoing flows, and use $c_q^{\a}$ to denote the concentration of gas component $\a$ at the node after mixing, which may be different from $c_q^{\alpha,s}$ in the case that $F_q^d=0$ and $F_q^s>0$, i.e., an injection at node $q$.  The condition on continuity of concentrations is then given as
\begin{align}\label{eq:continuity_of_d}
    \underline{c}_k^{\a} = c_q^{\a}, \qquad  \text{for all}\  k\in\partial^-q\  \text{and}\ \text{all}\  \a,
\end{align}
that is, the concentrations at the boundaries of pipes $k\in\partial^-q$ outgoing from a node $q\in\mathcal V$ are equal to the corresponding nodal concentrations after mixing of incoming flows.  The boundary conditions at a node $q$ that specify gas component mass fractions depend on whether there is an injection at that node.  If there is an \emph{outflow} from node $q$ such that $F_q^d\geq0$ and $F_q^s=0$, then the outflow concentrations are equal to the the nodal concentrations $c_q^{\alpha}$.  If there is an \textit{injection} at node $q$ such that $F_q^s\geq0$ and $F_q^d=0$, then the injection mass fractions must be specified as (possibly time-varying) parameters $c^\a_{q,s}$.  When stating an IBVP for pipeline simulation, the pressure is typically defined for at least one node of the graph, which is referred to as a ``slack node''.  This is in contrast to ``flow nodes'' where the flow in or out of the network is specified, as denoted by the flows $F_q^d$ (or $F_q^s$) withdrawn from (or injected into) the network as described above.   A well posed IBVP requires that either the pressure $p_q$ or flow $F_q^d$ (or $F_q^s$) is specified for each node $q\in\mathcal V$, and if  $F_q^s \geq 0$ then the mass fractions $c_q^{\alpha,s}$ must be specified for $n-1$ values of $\alpha$ where $n$ is the total number of gas constituents.  The pressure and flow boundary conditions are given by
\begin{subequations}
    \begin{align}\label{eq:pressure_and_flow_BC}
    p_j(t) = \text{ required parameter}, \,\,\, j\in\mathcal{V} \text{ is a slack node}, \\ F_q(t) = \text{ required parameter}, \,\,\, q \in\mathcal{V} \text{ is a flow node}.
\end{align}
\end{subequations}
The concentration boundary condition type and definition is then determined by 
\begin{align}\label{eq:continuity_of_d_outflow}
    c_q^{\a,s}(t) = \text{required parameter}, & \text{ if $F_q^s\geq 0$ and thus $F_q(t)\leq0$}.
\end{align}
Furthermore, we include the action of gas compressors in our model as nodal objects that boost pressure between the node and the start of a pipe.  If there is a compressor at the start of a pipe $k=(i,j)\in\mathcal{E}$ located at node $i\in\mathcal{V}$, then its action is defined by
\begin{align}\label{eq:compressor_action}
    \underline{p}_k(t) = \mu_k(t) p_i(t), 
\end{align}
where $\mu_k$ denotes the time-dependent compression ratio.  For the majority of pipes that do not have gas compressors at the start, we suppose that $\mu_k\equiv 1$.  Finally, in order to specify a well-posed IBVP for pipeline simulation, we also require initial conditions of flow and densities on each pipe in the graph.  These conditions take the form
\begin{align} \label{eq:initial_conditions}
    d_k^{\a}(0,x) = d_{k,0}^{\alpha}(x), \,\,\, \text{for all}\  k\in \mathcal E, \qquad
    \phi(0,x) = \phi_{k,0}(x), \,\,\, \text{for all}\  k\in \mathcal E.
\end{align}
It is critical that the initial conditions in equation \eqref{eq:initial_conditions} are compatible in the sense that boundary condition requirements in equations \eqref{eq:continuity_of_d}, \eqref{eq:continuity_of_d_outflow}, and \eqref{eq:compressor_action} are satisfied at the initial time $t=0$.
The equations \eqref{eq:multigaspde3}-\eqref{eq:initial_conditions} constitute a well-posed IBVP for heterogeneous gas flow over a pipeline network.  In the following section, we present a discretization scheme to be used for numerical solution of this type of simulation problem.

\section{Staggered Grid Discretization Scheme}  \label{sec:discretization}

We now present a numerical scheme for finite difference discretization of the equations \eqref{eq:multigaspde1a}-\eqref{eq:multigaspde1b}, which is adapted from the scheme for a single pipe that is described in a previous study \cite{gyrya2019explicit}.  We subsequently adapt the scheme for equations \eqref{eq:multigaspde2a}-\eqref{eq:multigaspde2b} for non-ideal gas flows and then to the treatment of nodal boundary conditions at junctions of multiple pipes in order to facilitate generalization to the setting of heterogeneous gas mixing as in the network IBVP in equations \eqref{eq:multigaspde3}-\eqref{eq:initial_conditions}. 

\subsection{Explicit Staggered-Grid Discretization for a Single Pipe} \label{sec:staggered_onepipe}

We consider a staggered finite difference (FD) discretization in which space grid points are indexed by two indices $i$ and $j$ and time layers are indexed by two additional indices $n$ and $m$.  Each density function $d^{\a}$ is discretized at the points $(t_n,x_i)$ and the total flow $\phi$ is discretized at the points $(t_m,x_j)$, where the first and second coordinates are reserved for time and space.  The staggered grids are defined as 
\begin{subequations}
\begin{align} 
    x_i = i\dx_k, &\quad i=\tfrac{1}{2},\tfrac{3}{2}, \dots,N_k-\tfrac{1}{2} 
    \quad \text{and} \quad \label{eq:grid_points_xi} \\
    x_j = j\dx_k, &\quad j=0,1, \dots,N_k,
    \quad \dx_k=L_k/N_k \label{eq:grid_points_xj}
\end{align}
\end{subequations}
and 
\begin{subequations}
\begin{align}
    t_n = n\dt, &\quad n= 0,1, \dots,M,
    \quad \text{and} \quad \label{eq:grid_points_ti} \\
    t_m = m\dt, &\quad m=\tfrac{1}{2},\tfrac{3}{2}, \dots,M+\tfrac{1}{2}, \label{eq:grid_points_tj}
\end{align}
\end{subequations}
where the number $N_k$ of discretization points may depend on edge $k$.  We henceforth denote the space discretization by $\dx$, with the understanding that it is in general pipe-dependent. The discrete representation of the \emph{conservation of mass} equation \eqref{eq:multigaspde1a}
for gas $\a$ at a point $(t_m,x_i)$, with staggered points $m=n+\tfrac{1}{2}$ and $i=j-\tfrac{1}{2}$, takes the form
\begin{align}\label{eq:conservation_of_mass_discrete}
	\tfrac{1}{\dt}(d_{n+1,i}^{\a} - d_{n,i}^{\a}) 
	+ 
	\tfrac{1}{\dx}
	\left( F^{\a}_{m,j} - F^{\a}_{m,j-1} \right)
	= 
	0,
\end{align}	
where the per-area component mass flux $\phi^{\a}_{m,j}$ through the pipe is defined as 
\begin{align}\label{eq:mass_flux_upwinding}
    \phi^{\a}_{m,j} 
    = 
    \left\{
    \begin{array}{ll}
      \phi_{m,j} \frac{d^{\a}_{n,i}}{d_{n,i}} \quad & \text{if }\ \phi_{m,j}\geq 0,  \\ \\
      \phi_{m,j} \frac{d^{\a}_{n,i+1}}{d_{n,i+1}} \quad & \text{if }\ \phi_{m,j}<0.  
    \end{array}
    \right.
\end{align}
The discrete representation of the \emph{conservation of momentum} equation \eqref{eq:multigaspde1b}
at a point $(t_{n+1},x_j)$, where $m=n+\tfrac{1}{2}$ and $i=j-\tfrac{1}{2}$, 
takes the form
\begin{align}\label{eq:conservation_of_momentum_discrete}
	\tfrac{1}{\dt}\left(\phi_{m+1,j} - \phi_{m,j}\right) 
	+ 
	\tfrac{1}{\dx}
	\left(p_{n+1,i+1} - p_{n+1,i} \right) 
	= 
	-
	f(\phi_{n+1,j},d_{n+1,j}),
\end{align}
where the pressure as defined in equation \eqref{eq:explicit_pressure} is
\begin{align*}
    p_{n+1,i+1} = p(d^{\a}_{n+1,i+1}), 
\end{align*}
and the friction term takes the form
\begin{align}\label{eq:friction_discrete}
    f(\phi_{n+1,j},d_{n+1,j})
    =
    \frac{\lambda}{2D}\cdot
    \frac{\phi_{m+1,j}|\phi_{m+1,j}| + \phi_{m,j}|\phi_{m,j}|}{ (d_{n+1,i+1} + d_{n+1,i})}.
\end{align}
Due to the form of the friction term \eqref{eq:friction_discrete},
the equation \eqref{eq:conservation_of_momentum_discrete} is a quadratic equation in $\phi_{m+1,j}$,
\begin{align}\label{eq:quadratic_equation}
    a\cdot \sign(\phi_{m+1,j})\cdot (\phi_{m+1,j})^2 + (\phi_{m +1,j}) - c = 0,
\end{align}
with coefficients $a$ and $c$ given by
\begin{align*}
    a &= \dt \cdot \frac{\lambda L}{2D} \cdot \frac{1}{d_{n+1,i+1} + d_{n+1,i}},\\
    c &= \phi_{m,j} 
         - \frac{\dt}{\dx} \left(p_{n+1,i+1} - p_{n+1,i} \right) 
         - \dt \cdot \frac{\lambda L}{2D}\cdot \frac{\sign(\phi_{m,j})\cdot (\phi_{m,j})^2}{d_{n+1,i+1} + d_{n+1,i}}.
\end{align*}
The relevant solution of the quadratic equation \eqref{eq:quadratic_equation} is
\begin{equation}\label{eq:explicit_solution_for_flux_phi}
    \phi_{m+1,j}
    %=F^{-1}(y) 
    = \sign(c)\frac{-1 + \sqrt{1+4a|c|}}{2a}.
\end{equation}

\subsection{Nodal Boundary Conditions} \label{sec:Nodal_BC}

The nodes of the gas network graph represent the junctions that in general connect multiple pipes. Recall that as defined in Section \ref{sec:network}, we denote the set of pipes oriented as \emph{entering} node $q$ by $\partial^+q$ and the set of pipes oriented as \emph{leaving} node $q$ as $\partial^-q$. Let us then use $\partial q = \partial^-q \cup \partial^+q$ to denote all pipes attached to the node $q$.
We will use the superscript $k$ to enumerate the pipes adjacent to the node and to indicate that particular quantity belongs to the pipe $k$.

In our notation, a discretized variable on an edge $k\in\mathcal{E}$ is denoted and indexed by, e.g., $d_{m,j}^{\alpha,k}$ for $d^\alpha$, where  $\alpha$ indicates the gas component, $k$ indicates the edge, $m$ indicates time, and $j$ indicates position on the edge.  We may write, e.g., $d_{m,j=0}^{\alpha,k}$ or $d_{m,i=1/2}^{\alpha,k}$ to refer to the edge endpoint value on the main or staggered grid, respectively, and we may write, e.g., $d_{m-1/2,i=1/2}^{\alpha,k}$ to indicate position on the staggered time grid.  The $\alpha$ index is not used when referring to the total density of the mixed gas.  Similarly, a discrete-time variable at time $t_n$ is denoted and indexed by, e.g., $p_n^{q}$ for pressure and $d_n^{\alpha,q}$ for component density, where $q\in\mathcal V$ denotes the node, $\alpha$ indicates the gas component, and $n$ indicates time.  We may write $p_{m-1/2}^{q}$ to indicate position on the staggered time grid.  We similarly write $\mu_m^k=\mu_k(t_m)$ as the discrete time version of the compression ratio of a compressor at the start of a pipe $k$.
We will write total injection flow as $F_m^{q}$, component injection flow as $F_m^{\alpha,q}$, total injection density as $d_{m}^{q}$, and component injection density as $d_{m}^{\alpha,q}$.

To formulate the discrete version of the nodal boundary conditions, we start with the conservation of mass conditions \eqref{eq:flow_balance_general_form} for each of the gas components at time $t_m$, which are given for a node $q$ as
%Mass balance is modeled by classical Kirchhoff-type constraints
\begin{equation}\label{eq:disc_flow_balance_at_node}
    \sum_{k\in\partial^+q} S_k \phi_{m,j=N_k}^{\a,k}-\sum_{k\in\partial^-q} S_k \phi_{m,j=0}^{\a,k} = F_m^{\a,q,d} - F_m^{\a,q,s}, \quad \forall \,\alpha\in\{1,\ldots,n\}, \, \forall \, q\in\mathcal{V},
\end{equation}
where $S_k=\pi (D_k/2)^2$ is the cross-sectional area of pipe $k$,
$\phi_{m,j=0}^{\a,k}$ and $\phi_{m,j=N_k}^{\a,k}$ denote the mass flux of gas $\a$ per unit area
at the points $x_{0}=0$ and $x_{N_k}=L_k$  of the pipe $k$,  
as defined in \eqref{eq:mass_flux_upwinding}.  Here $F_m^{\a,q,d}$ and $F_m^{\a,q,s}$ are the rates of withdrawal or injection of gas $\a$ from (into) the pipeline system at node $q$ at time $t_m$, which correspond to the terms $c^\alpha(t_m)F_q^d(t_m)$ and $c^{\alpha,s}(t_m)F_q^s(t_m)$ in equation \eqref{eq:flow_balance_general_form}, respectively.
For simplicity of exposition in describing nodal conditions at a node $q$, we will henceforth suppose that local enumeration of the space grid points associated with each pipe $k\in\partial q$ begins at the node being examined, i.e. $j=0$ corresponds to the start of the pipe which is adjacent to that node. We may then re-write the flow balance laws \eqref{eq:disc_flow_balance_at_node} as
\begin{equation}\label{eq:disc_flow_balance_at_node_simple}
    \sum_{k\in\partial q} S_k \phi_{m,j=0}^{\a,k} + F_m^{\a,q,d} - F_m^{\a,q,s} = 0, \quad \forall \,\alpha\in\{1,\ldots,n\}, \, \forall \, q\in\mathcal{V}.
\end{equation}
Observe that the individual gas component withdrawals $F_m^{\a,q,d}$ and injections $F_m^{\a,q,s}$ are specified as parameters in the conditions \eqref{eq:disc_flow_balance_at_node} with the same type of complementarity conditions as for \eqref{eq:flow_balance_general_form}.   That is, $F_m^{\a,q,s}\geq 0$ and $F_m^{\a,q,d}=0$ in the case of gas injection, $F_m^{\a,q,d}\geq0$ and $F_m^{\a,q,s}=0$ if gas is withdrawn.  In the latter case, the withdrawal mass flow rates of individual components depend on the gas composition at the node $q$.  We denote positive mixed gas withdrawal at node $q$ by $F_m^{q,d}\equiv F_{q}^{d}(t_m)$,  and the withdrawal flow of gas component $\alpha$ computed as 
\begin{align}\label{eq:F_alpham}
    F_m^{\a,q,d} = \frac{d^{\a,q}_{m-1/2}}{d_{m-1/2}^{q}}\cdot F_m^{q,d},
    \qquad \text{when } F_m^{q,d}\geq 0,
\end{align}
where $d^{\a,q,d}_{m-1/2}$ are the densities of the individual gas components injected at the node and $d^{q,d}_{m-1/2}=\sum_\a d^{\a,q,d}_{m-1/2}$ is the derived total gas density of gas injected into the node.  We evaluate densities on the staggered time grid at $m-1/2$, a half step prior to time $m$ on which flows are evaluated. We describe computation of nodal component densities $d^{\a,q,d}_{m-1/2}$ subsequently.
\begin{figure}[ht!]
    \centering
    \includegraphics[width=.85\textwidth]{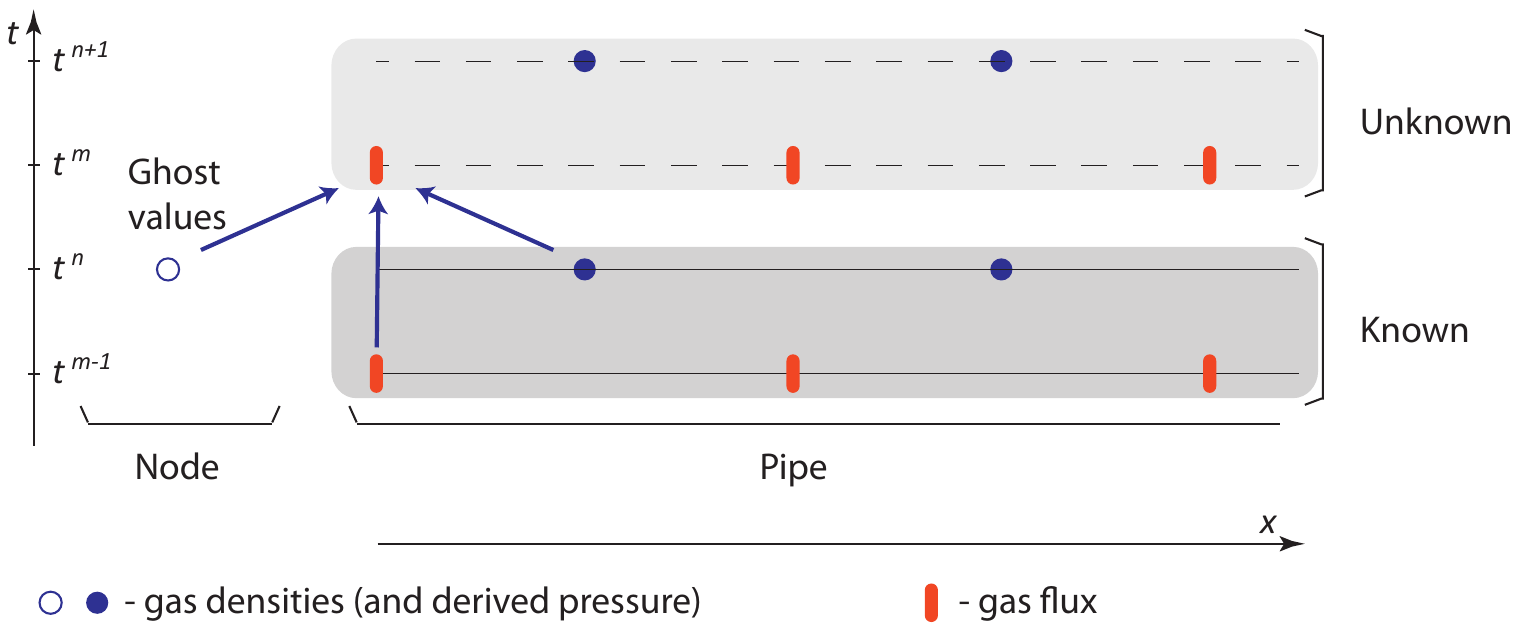}
    \caption{Illustration of the dependencies in the double step of time integration near the nodal boundary.}
    \label{fig:nodal BC}
\end{figure}

Observe in the case when $\phi_{m,j=0}^{k}>0$ that by definition \eqref{eq:mass_flux_upwinding} computation of
the flux component $F_{m,j=0}^{\a,k}$ of gas $\alpha$ requires knowledge of $\phi_{m,j=0}^{k}$ and $d_{m-1/2,i=-1/2}^{\a,k}$.
We will associate the values $d_{m-1/2,i=-1/2}^{\a,k}$ for all $k$ attached to a node $q$ and require them to be equivalent for all $k\in\partial q$ to the auxiliary nodal gas constituent densities $d_{m-1/2}^{\a,q}$ associated with the node $q$, so that 
\begin{align} \label{eq:node_edge_density_grid_compatibility}
    d_{m-1/2}^{\a,q} = d_{m-1/2,i=-1/2}^{\a,k},
    \qquad \text{for all}\  k\in\partial q.
\end{align}
The dependencies are illustrated in a schematic in Fig.~\ref{fig:nodal BC}.  The application of the nodal boundary conditions thus
requires determining $d^{\a,q}_{m-1/2}$, 
associated with the node $q$ at the time $t_{m-1/2}$, such that the flows $\phi_{m,j=0}^{k}$
computed from the discrete conservation of momentum equation \eqref{eq:conservation_of_momentum_discrete} 
satisfy the conservation of mass equations \eqref{eq:disc_flow_balance_at_node}. Because the resulting system of equations becomes increasingly complicated as the number of constituent gases increases, we reformulate the boundary conditions so that they are defined entirely in terms of nodal pressure $p_n^{q}$, at node $q$ and at time $t_n$, as follows.

The discrete version of the \emph{conservation of momentum} equation \eqref{eq:multigaspde1b} at the start of a pipe $k$ outgoing from node $q$ with a compressor is defined at a point $(t_{n},x_{j=0})$, with staggered points defined by $m=n+\tfrac{1}{2}$ and $i=j-\tfrac{1}{2}$, takes the form
\begin{align}\label{eq:conservation_of_momentum_discrete_alt}
	\frac{1}{\dt}\left(\phi_{m,j=0}^{k} - \phi_{m-1,j=0}^{k}\right) 
	+ 
	\frac{1}{\dx}
	\left(p^{k}_{n,i=1/2} - \mu_n^k p_n^{q} \right) 
	= 
	-
	f(\phi_{n,j=0}^{k},d_{n,j=0}^{k}),
\end{align}
where the pressure is
\begin{align*}
	p^{k}_{n,i=1/2} = p(\Bd^{\a,k}_{n,i=1/2}), 
\end{align*}
as defined in \eqref{eq:explicit_pressure}, and where $\mu_m^k \ge 1$ is a compression ratio that multiplies the suction pressure $p_n^{q}$ at node $q$ to yield the discharge pressure of the compressor at the start of the pipe, and the friction term is
\begin{align}\label{eq:friction_discrete_alt}
	f(\phi_{n,j=0}^{k},d_{n,j=0}^{k})
	=
	\frac{\lambda_k}{2D_k}\cdot
	\frac{\phi_{n,j=0}^{k}|\phi_{n,j=0}^{k}|}{d_{n,j=0}^{k}}.
\end{align}
Both $\phi$ and $d$ must be evaluated between staggered grid points where their values are in practice not available. Therefore, we approximate these values by upwinding to obtain
\begin{align}\label{eq:friction_term_upwinding}
	f(\phi_{n,j=0}^{k},d_{n,j=0}^{k}) 
	= 
	\frac{\lambda_k}{2D_k}\cdot
	\frac{\phi_{m-1,j=0}^{k}|\phi_{m-1,j=0}^{k}|}{d_{n,i=1/2}^{k}}.  
\end{align}
It follows that substituting \eqref{eq:conservation_of_momentum_discrete_alt} into the summation of \eqref{eq:disc_flow_balance_at_node_simple} over $\alpha$ and reformulating yields an explicit expression for the nodal pressure $p_n^{q}$ at a node $q$ at time $t_n$, given by 
\begin{align}\label{eq:flow_balance_scalar_coef1}
	& p_n^{q} \!=\! \frac{\dx}{\dt}\!\cdot\! \frac{ \left( F^{q,s}_m-F^{q,d}_m - \sum_{k\in\partial q} S_k  \left(\phi_{m-1,j=0}^{k} - \frac{\dt \lambda_k}{2D_k} \cdot \frac{\phi_{m-1,j=0}^{k}|\phi_{m-1,j=0}^{k}|}{\sum_\a d_{m - 1/2,i = 1/2}^{\a,k}} - \frac{\dt}{\dx} p_{n,i=1/2}^{k} \right)\right)}{\sum_{k\in\partial q} S_k\mu_n^k},
\end{align}	
where $\mu_n^k\equiv 1$ for all pipes $k$ with no compressors.  Detailed calculations for both ideal and non-ideal gas modeling can be found in Appendix \ref{sec:append_nodal_bc}. 

In addition, we can obtain values for mass flux and densities at the boundary of a pipe expressed in terms of nodal pressure $p_n^{q}$ at the boundary node $q$, which we write as $\phi^{k}_{m,j=0}\left(p_n^{q}\right)$ and $d^{\a,k}_{m + 1/2,i = 1/2}\left(p_n^{q}\right)$. We can do this for both ideal and non-ideal gas models by using the results of Section \ref{sec:EOS} to make appropriate substitutions into equation \eqref{eq:conservation_of_momentum_discrete_alt} for pressure variables as functions of component densities.
For an ideal gas, the boundary mass flux into pipe $k$ from node $q$ with pressure $p_n^{q}$ at time $t_n$ is given by
\begin{equation}\label{eq:boundary_flow}
	\phi^{k}_{m,j=0}\left(p_n^{q}\right)
	\!= \!\left(\!\phi_{m-1,j=0}^{k}
	\!-\!
	\frac{\dt \lambda_k}{2D_k}
	\frac{\phi_{m-1,j=0}^{k}|\phi_{m-1,j=0}^{k}|}{\sum_\a d^{\a,k}_{m - 1/2,i = 1/2}} \!-\! \frac{\dt}{\dx} \! \sum_\a d^{\a,k}_{m - 1/2,i = 1/2} R_\a T \! \right) \!+\! \frac{\dt}{\dx} \mu_n^kp_n^{q}.
\end{equation}
The above expression can then in turn be used together with the pressure-density relation $d_n^{\alpha,q}(p_n^q)=\gamma_n^{\alpha,q}\pi_\a^{-1}(p_n^q)$ to obtain the update formula boundary values of densities of the gas components for flow leaving node $q$ into the pipe $k$:
\begin{align}\label{eq:conservation_of_mass_discrete_ideal}
	d^{\a,k}_{m + 1/2,i = 1/2}\left(p_n^{q}\right) \!=\! d^{\a,k}_{m - 1/2,i = 1/2} 
	\!-\! 
	\frac{\dt}{\dx}
	\!\left( \phi^{k}_{m,j=1} \frac{ d^{\a,k}_{m - 1/2,i = 1/2}}{\sum_\a d^{\a,k}_{m - 1/2,i = 1/2}} \!-\! \phi^{k}_{m,j=0}\left(p_n^{q}\right) \frac{d_n^{\alpha,q}\left(p_n^{q}\right)}{\sum_\a d_n^{\alpha,q}\left(p_n^{q}\right)} \right)\!\!.
\end{align}
The above equations \eqref{eq:boundary_flow} and \eqref{eq:conservation_of_mass_discrete_ideal} can be modified to accommodate non-ideal gas modeling by using the expression for pressure in equation \eqref{eq:explicit_pressure}.  This results in an expression for boundary mass flux into pipe $k$ from node $q$ with pressure $p_n^{q}$ at time $t_n$ is given by
\begin{align}\label{eq:boundary_flow_nonideal}
	&\phi^{k}_{m,j=0}\left(p_n^{q}\right)
	\!=\! \left(\phi_{m-1,j=0}^{k}
	\!-\!
	\frac{\dt \lambda_k}{2D_k}
	\frac{\phi_{m-1,j=0}^{k}|\phi_{m-1,j=0}^{k}|}{\sum_\a d^{\a,k}_{m - 1/2,i = 1/2}} \!-\! \frac{\dt}{\dx} \frac{\sum_\a d^{\a,k}_{m - 1/2,i = 1/2} R_\a T}{1 - \sum_\a d^{\a,k}_{m - 1/2,i = 1/2} R_\a T a_\a} \right) \!+\! \frac{\dt}{\dx} p_n^q.
\end{align}
The expression for updating the boundary values of densities of the gas components for flow leaving node $q$ into the pipe $k$ will be the same as equation \eqref{eq:conservation_of_mass_discrete_ideal}:
\begin{equation}\label{eq:conservation_of_mass_discrete_nonideal}
	d^{\a,k}_{m + 1/2,i = 1/2}\left(p_n^{q}\right) \!=\! d^{\a,k}_{m - 1/2,i = 1/2} 
	\!-\! 
	\frac{\dt}{\dx}
	\!\left( \phi^{k}_{m,j=1} \frac{ d^{\a,k}_{m - 1/2,i = 1/2}}{\sum_\a d^{\a,k}_{m - 1/2,i = 1/2}} \!-\! \phi^{k}_{m,j=0}\left(p_n^{q}\right) \frac{d_n^{\alpha,q}\left(p_n^{q}\right)}{\sum_\a d_n^{\alpha,q}\left(p_n^{q}\right)} \right)\!\!. \!
\end{equation}
The above explicit discretization scheme can be applied to solve the IBVP for heterogeneous gas mixture flow in a pipeline network with compressors defined by equations \eqref{eq:multigaspde3}-\eqref{eq:initial_conditions}. A more detailed derivation of the above results is provided in Appendix \ref{sec:append_nodal_bc}.

\section{Nodal Monitoring Policies} \label{sec:monitoring}

A recent study finds that blending more than five percent hydrogen by volume into natural gas pipelines results in a greater likelihood of pipeline leaks \cite{cpuc2022impacts}. Blends of natural gas with greater than 5\% hydrogen by volume could require modifications of appliances such as stoves and water heaters to avoid leaks and equipment malfunction. Blending more than 20\% hydrogen by volume presents a higher likelihood of permeating plastic pipes, which can increase the risk of gas ignition outside the pipeline. Further problems with the transportation of hydrogen using industrial pipelines that were originally designed to transport natural gas have been known for some time \cite{melaina2013blending}.  Here we propose two real-time monitoring and corrective flow control policies.  The first can be used to ensure that hydrogen injections into a natural gas pipeline system do not cause mass or volumetric fractions to exceed allowable limits, and the second can be used to ensure that hydrogen injections do not cause pressures to fall below minimum requirements.

\subsection{Input Flow Monitoring for Maximum Concentration Limit}  \label{sec:input_flow_monitoring_concentration}
After being equipped with an explicit numerical method for simulation of pipeline flows with heterogeneous gas injections, we can develop a local nodal monitoring policy to guarantee that upper bounds on hydrogen mass fraction are satisfied in real time.  Consider the following scenario.  Suppose that there is a gas mixture flowing at a constant rate into a given node $q$ with known individual gas concentrations.  Let us suppose that the composition of this injection flow is different from the composition of gases flowing into the node from incoming network pipes $k\in\partial^+q$, and the nodal injection thus alters the composition of the flow leaving the node to outgoing pipes $k\in\partial^-q$.  Suppose that at some point, the mass fraction $c_m^{\alpha,q}$ at node $q$ of a constituent gas $\alpha$ approaches a safety limit $c_{\max}^{\alpha,q}$.  In particular, we consider a scenario in which the volumetric concentration of hydrogen in a mixture with natural gas must not exceed, e.g., 10\%, and the injection at node $q$ is hydrogen gas.  The monitoring and control policy goal is to control the gas injection at node $q$ in order to ensure that the safety requirements are satisfied.

By substituting the expression for $p_n^q$ in equation \eqref{eq:flow_balance_scalar_coef1} into equation \eqref{eq:boundary_flow}, it is possible to express the boundary flow $\phi^{k}_{m,j=0}$ as a function of the injection flow $F_m^{q,s}$, %as $\phi^{k,m}_{j=0}\left(F_m^{q,s}\right)$, 
i.e.,
\begin{align}\label{eq:flow_input_flow_dependence_2}
         \phi^{k}_{m,j=0}(F_n^{q,s}) = \frac{\mu_n^k F_n^{q,s}}{\sum_r \mu_n^r S_r} + \Theta_m^{q,k} + \mu_n^k \!\cdot\! \Upsilon_m^q.
\end{align}
where we assume that withdrawal $F_m^{q,d}$ is zero, and $\Theta_m^{q,k}$ and $\Upsilon_m^q$ are given by
\begin{align} \label{eq:outgoing_flow_term_thetamk}
    & \Theta_m^{q,k} = \left(\!\phi_{m-1,j=0}^{k}
	-
	\frac{\dt \lambda_k}{2D_k}
	\frac{\phi_{m-1,j=0}^{k}|\phi_{m-1,j=0}^{k}|}{\sum_\a d^{\a,k}_{m - 1/2,i = 1/2}} \!-\! \frac{\dt}{\dx} p_{n,i=1/2}^{k} \!\right), \quad \forall k \in \partial q,\\
	& \Upsilon_m^q = \frac{ \sum_{k\in\partial q} S_k  \left( \frac{\dt \lambda_k}{2D_k} \cdot \frac{\phi_{m-1,j=0}^{k}|\phi_{m-1,j=0}^{k}|}{\sum_\a d_{m - 1/2,i = 1/2}^{\a,k}} - \frac{\dt}{\dx} p_{n,i=1/2}^{k} - \phi_{m-1,j=0}^{k}\right)}{\sum_{k\in\partial q} S_k\mu_n^k}. \label{eq:outgoing_flow_term_upsilonmq}
\end{align}
A more detailed derivation of equations \eqref{eq:flow_input_flow_dependence_2}-\eqref{eq:outgoing_flow_term_upsilonmq} and further analysis can be found in Appendix \ref{sec:append_monitoring}.  Following the subsequent derivation in Appendix \ref{sec:append_monitoring}, and recalling that mass fraction at a node $q$ and on a pipe $k$ is defined as
\begin{equation} \label{eq:concentration_as_density_fraction}
    c_m^{\a,q}=\frac{d_m^{\alpha,q}}{\sum_{\beta} d_m^{\beta,q}}, \,\, \forall q\in\mathcal{V}, \quad c_{m,j}^{\a,k}=\frac{d_{m,j}^{\alpha,k}}{\sum_{\beta} d_{m,j}^{\beta,k}}, \,\, \forall k\in\mathcal{E},
\end{equation} 
we obtain an expression in equation \eqref{eq:max_injection_at_mass_fraction_limit} for the maximum amount $F_{\max}^{q,s}$ of hydrogen injection that maintains the nodal hydrogen mass fraction just below the maximum allowable limit:
\begin{equation} \label{eq:max_injection_at_mass_fraction_limit}
F_{\max}^{q,s}(c^{\alpha,q}_{\max}) \!=\! \frac{c_{m-1/2, i = 1/2}^{\a,k} \sum_{k\in\partial_+r} S_k \left(\Theta_m^{q,k} + \mu_n^k \!\cdot\! \Upsilon_m^q \right) \!-\! c^{\alpha,q}_{\max} \sum_{k\in\partial_-r}  S_k \left( \Theta_m^{q,k} + \mu_n^k \!\cdot\! \Upsilon_m^q \right) }{  c^{\alpha,q}_{\max} \sum_{k\in\partial_-r}  S_k \left( \frac{\mu_n^k}{\sum_r \mu_n^r S_r} \right) - c_{m-1/2, i = 1/2}^{\a,k} \sum_{k\in\partial_+r} S_k \left( \frac{\mu_n^k}{\sum_r \mu_n^r S_r} \right) - c_m^{\alpha,s}}.
\end{equation}
The injection at node $q$ must be bounded above zero, i.e., in the case that the mass fraction of hydrogen flowing into the node is already at the allowable maximum $c_{\max}^{\alpha,q}$.  The injection would be at the pre-planned level $F_m^{q,s}$ when $0 < F_m^{q,s}< F_{\max}^{q,s}$.

\subsection{Output Flow Monitoring for Minimum Pressure Requirement} \label{sec:input_flow_monitoring_pressure}

Engineering limitations of pipeline systems require gas pressures to be maintained between minimum pressures needed for customers to withdraw the gas, and maximum limits that cannot be exceeded in order to maintain material integrity.  As long as the pressures at slack nodes and the discharge pressures of compressors is at or below the allowable maximum for each pipe, the decrease of pressure in the direction of flow that is the characteristic of the usual regime of pipeline operations will ensure that maximum pressure limits are not exceeded.  However, small amounts of hydrogen injection may result in significant pressure decrease along a pipe because of the increased wave speed of the resulting mixture.  This may potentially result in local violations of minimum pressure constraints at nodes with consumers.  Predicting allowable transient changes in hydrogen injections \emph{a priori} is prohibitively difficult, particularly for pipeline networks with complex topology that may include loops.

We mention here that by the monotonicity property of transient dissipative flows on graphs, minimum pressure on a pipe and thus at network nodes occurs at downstream locations where flows $F_m^{q,d}$ are withdrawn \cite{misra2020monotonicity}.  Thus, we observe that the decrease in pressure along a pipeline may be limited by limiting flow velocities, and this can be accomplished locally by limiting the withdrawal by relevant consumers.  Suppose then that the planned withdrawal $F_m^{q,d}$ at a node $q$ would result in pressure falling below the minimum required value $p_{\min}^q$ at that node. We address the issue of maintaining pressures above required location-dependent minimum values $p_{\min}^q$ by deriving a formula for the maximum withdrawal $F_{\max}^{q,d}$, which is less than the planned withdrawal $F_m^{q,d}$, that will maintain the nodal pressure at exactly $p_{\min}^q$ in that scenario.  The result is entirely local and intended for real-time response, and is not intended to be predictive or optimal.  The latter type of analysis requires optimization algorithms \cite{baker2023optimal}.

In order to derive the real-time monitoring policy for maintaining minimum pipeline pressures, recall that when $q$ is a withdrawal node and thus $F_m^{q,s}\equiv 0$ the equation \eqref{eq:flow_balance_scalar_coef1} for $p_n^q$ appears as
\begin{align} \label{eq:flow_balance_scalar_coef1_withdrawal}
	& p_n^{q} \!=\! \frac{\dx}{\dt}\!\cdot\! \frac{ \left(-F^{q,d}_m - \sum_{k\in\partial q} S_k  \left(\phi_{m-1,j=0}^{k} - \frac{\dt \lambda}{2D} \cdot \frac{\phi_{m-1,j=0}^{k}|\phi_{m-1,j=0}^{k}|}{\sum_\a d_{m - 1/2,i = 1/2}^{\a,k}} - \frac{\dt}{\dx} p_{n,i=1/2}^{k} \right)\right)}{\sum_{k\in\partial q} S_k\mu_n^k}.
\end{align}	
In the case that the planned withdrawal $F^{q,d}_m$ would lead to a value of $p_n^{q}$ that is less than the allowable minimum $p_{\min}^q$, the withdrawal can be curtailed to a maximum value obtained by re-arranging equation \eqref{eq:flow_balance_scalar_coef1_withdrawal} as
\begin{align}\label{eq:max_withdrawal_to_maintain_pmin}
	\!\!&  F_{\max}^{q,d} \!=\!  \!\sum_{k\in\partial q} \! S_k  \left(\frac{\dt}{\dx} p_{n,i=1/2}^{k} \!+\! \frac{\dt \lambda}{2D} \cdot \frac{\phi_{m-1,j=0}^{k}|\phi_{m-1,j=0}^{k}|}{\sum_\a d_{m - 1/2,i = 1/2}^{\a,k}} \!-\! \phi_{m-1,j=0}^{k} \right)\! -\frac{\dx}{\dt}\!\sum_{k\in\partial q} S_k\mu_n^k  p_{\min}^{q}.
\end{align}
By inspecting equation \eqref{eq:flow_balance_scalar_coef1_withdrawal}, we find that if a compressor is present at the start of a pipe $k$ directed outward from node $q$, then decreasing its compression ratio $\mu_n^k$ would be another way of maintaining the nodal pressure $p_n^q$ above the minimum required level.  We do not examine such recourse policies for two reasons.  First, changing compressor settings will significantly alter the downstream pressure dynamics thus causing cascading issues.  Second, compressor stations are in practice present at relatively few locations in a pipeline network that do not coincide with locations of major gas consumption, such as city gates.

\section{Computational Studies} \label{sec:computational}

We present a collection of numerical simulations to address open questions about modeling needs in the transport of gas mixtures, to demonstrate the functionality and generality of our numerical scheme, to cross-validate a benchmark simulation with respect to the results of a previous study, and to demonstrate the use of the nodal monitoring policy.

\subsection{Single Pipe Simulations} \label{sec:onepipe_example}

We first formulate an IBVP for gas flow through a single pipe, with dynamics, initial conditions, parameters, and boundary conditions given below.  Here $a_1$ and $a_2$ are the wave speeds of hydrogen and natural gas, respectively, and $\epsilon_1$ and $\epsilon_2$ are the respective diffusion coefficients.  The initial conditions represent steady flow of natural gas, and the mass fraction of the injection gas is given by $c(t)$, where $c(0)=0$.  This IBVP was examined in previous studies \cite{gyrya2019explicit,herty2010new} for homogenous gas only, and those results are used for benchmarking.

~

\begin{minipage}{0.45\textwidth}
    \text{Dynamics}: 
    \begin{equation*}
        \begin{cases} \partial_t d_1 + \partial_x \Big( \frac{d_1 \phi }{d_1 + d_2}\Big) = \epsilon_1 \Delta d_1 ,\\\partial_t d_2 + \partial_x \Big( \frac{d_2 \phi }{d_1 + d_2}\Big) = \epsilon_2 \Delta d_2,\\  \partial_t \phi + \partial_x (d_1 R_1 T + d_2 R_2 T) = -\frac{\lambda}{2D}\frac{\phi|\phi|}{d_1 + d_2}\end{cases}
    \end{equation*}\\
\end{minipage}
\begin{minipage}{0.45\textwidth}
    \text{Initial Conditions}:
    \begin{equation*}
         \begin{cases} d_1(0,x) = 0 \\ d_2(0,x) = \sqrt{{s(0)}^2 - \frac{\lambda}{(R_2 T) D} \cdot \phi_0 |\phi_0|x} \\ \phi(0,x) = \phi_0 \end{cases}
    \end{equation*}\\
\end{minipage}

~

\begin{minipage}{0.4\textwidth}
\text{Parameters}:
\begin{equation*}
    \begin{cases} \lambda = 0.011\\ a_1 = \sqrt{R_1 T} = 1320 \quad m/s \quad\\ a_2 = \sqrt{R_2 T} = 377.9683 \quad m/s \quad  \\ D=0.5 \quad m\\ L = 100 \times 10^3 \quad m\\ T = 3600 \times 12 \quad s\\ \phi_0 = 289 \quad kg/m^2 /s\\ \rho_0 = 45.4990786148 \quad kg/m^3\\ s(t) = \rho_0 \cdot (1 + \frac{1}{10}\sin(6\pi t/T))\\ l(t) = \phi_0 \cdot (1 + \frac{1}{10}\sin(4\pi t/T))\end{cases}
\end{equation*}\\
\end{minipage}
\begin{minipage}{0.15\textwidth}
    \text{Variables}:
    \begin{equation*}
        \begin{cases} d_1 : [0, T] \times [0, L] \rightarrow \mathbb{R}_{+} \\d_2 : [0, T] \times [0, L] \rightarrow \mathbb{R}_{+} \\ \phi : [0, T] \times [0, L] \rightarrow \mathbb{R}\\ t \in [0,T], \quad x \in [0,L] \end{cases}
    \end{equation*}
    \text{Boundary Conditions}: 
    \begin{equation*}
        \begin{cases} d_1(t,0) = s(t)  \cdot (1 - c(t)) \\ d_2(t,0) = s(t) \cdot c(t)  \\ c(t) = 0.4t/T, ~~ \mbox{if}  ~~ t \le T/4, ~~ \mbox{else} ~~ c(t) = 0.1 \\ \phi(t,L) = l(t) \end{cases}
        %\begin{cases} d_1(t,0) = (t<T/4)\cdot(0.1 t/(T/4)) + (t \ge %T/4)(0.1) s(t)\\ d_2(t,0) = (t<T/4)(0.1 t/(T/4)) + (t \ge %T/4)(0.1) (1- s(t))\\ \phi(t,L) = l(t) \end{cases}
    \end{equation*}\\
\end{minipage}

\begin{figure}[ht!]
    \centering
    \vspace{-2ex}
	  \includegraphics[width=1\linewidth]{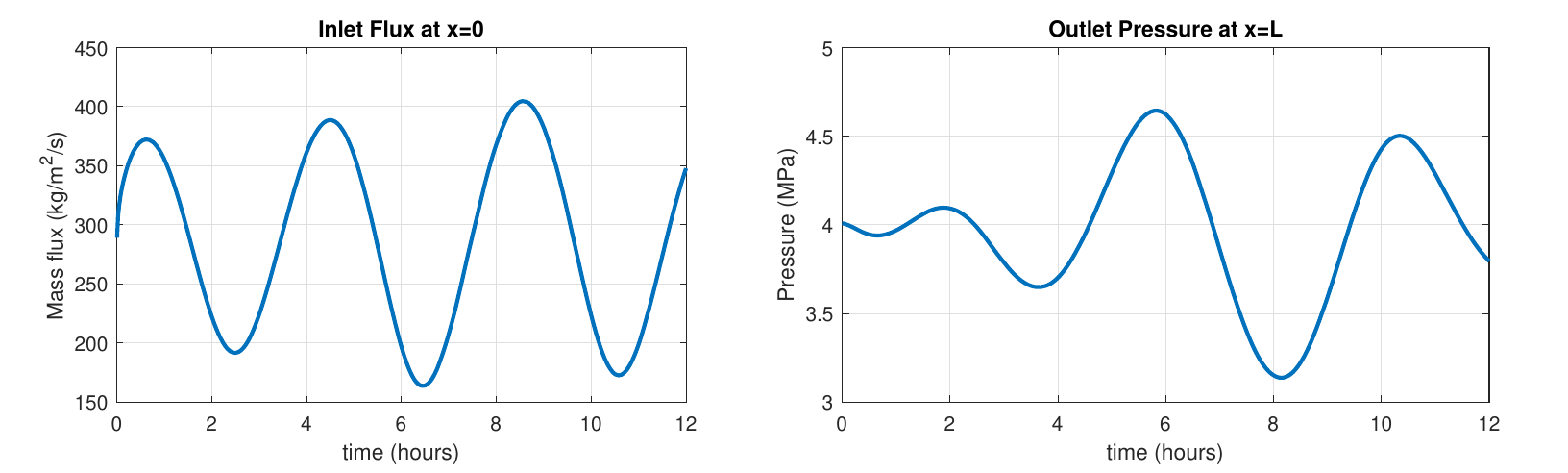}
       \vspace{-4ex}
	  \caption{Single pipe simulation when $c(t)\equiv 0$ so no hydrogen is injected. Left: inlet flux; Right: outlet pressure.}
   \label{fig:one_pipe_one_gas}
\end{figure}

The input pressure and outlet mass flux are given by $s(t)$ and $l(t)$ in the parameter definition above, and the input flux and outlet pressure are results of the simulation.  All simulations are implemented directly in MATLAB without relying on any packages or third-party solvers.  For the single pipe simulations, 200 space discretization points and 40000 time discretization points are used. The time required for each simulation is approximately $t=1.2$ seconds, although we note that our implementation is not optimized for rapid computation.  We first conduct the simulation with no hydrogen injection and no diffusion modeling, that is, with $c(t)\equiv 0$ and $\epsilon_1=\epsilon_2=0$, and the results are shown in Figure \ref{fig:one_pipe_one_gas}.
\begin{figure}[ht!]
    \centering
    \vspace{-2ex}
	  \includegraphics[width=\linewidth]{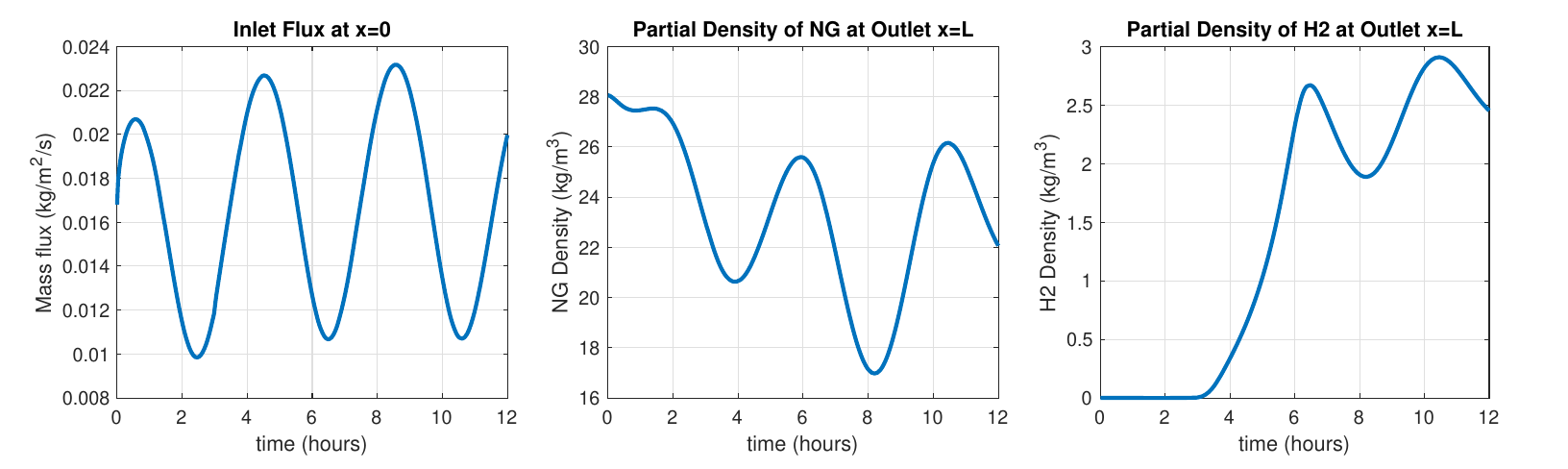}
   \vspace{-4ex}
	  \caption{Single pipe simulation with hydrogen blending according to $c(t)$ as defined above, up to 10\% of added hydrogen. Left: total inlet flux at $x=0$; Center: partial density of natural gas at outlet $x=L$; Right: partial density of hydrogen at outlet $x=L$.}
      \label{fig:one_pipe_mixing}
\end{figure}

\begin{figure}[ht!]
    \centering
        \vspace{-1ex}
	  \includegraphics[width=1\linewidth]{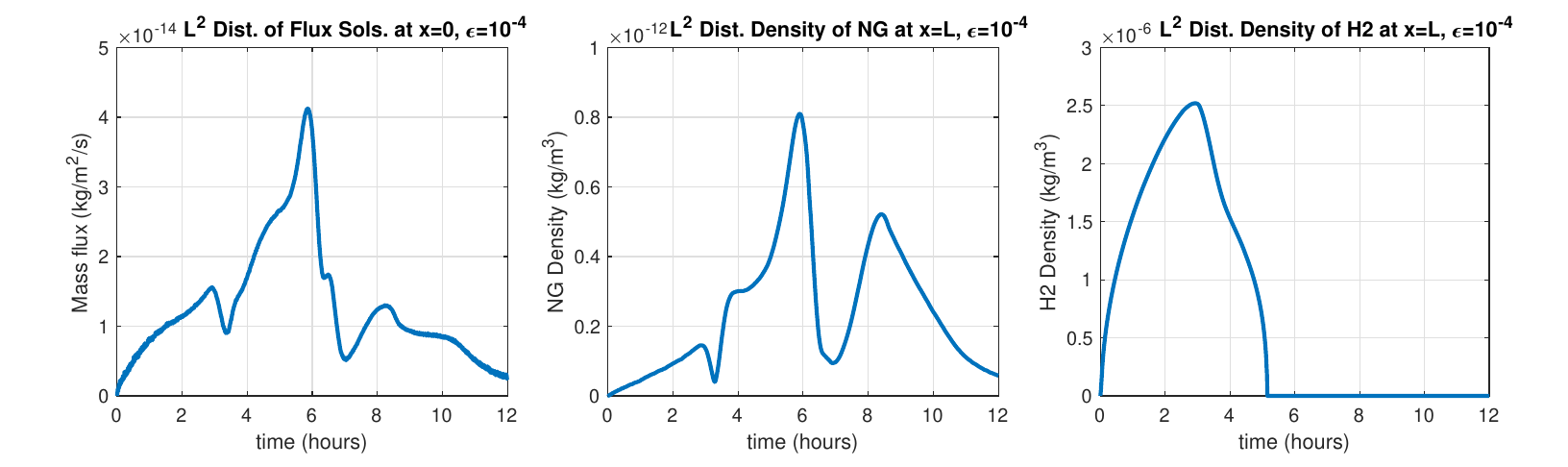}
       \vspace{-4ex}
	  \caption{Plots of $L^2$ differences in simulations with and without accounting for diffusion with  coefficients $\epsilon_1=\epsilon_2=0.0001$.  Left: difference in input flux; Center: natural gas densities; Right: hydrogen densities.}
       \label{fig:one_pipe_mixing_mismatch_00001}
\end{figure}

\begin{figure}[ht!]
    \centering
     \vspace{-1ex}
	  \includegraphics[width=1\linewidth]{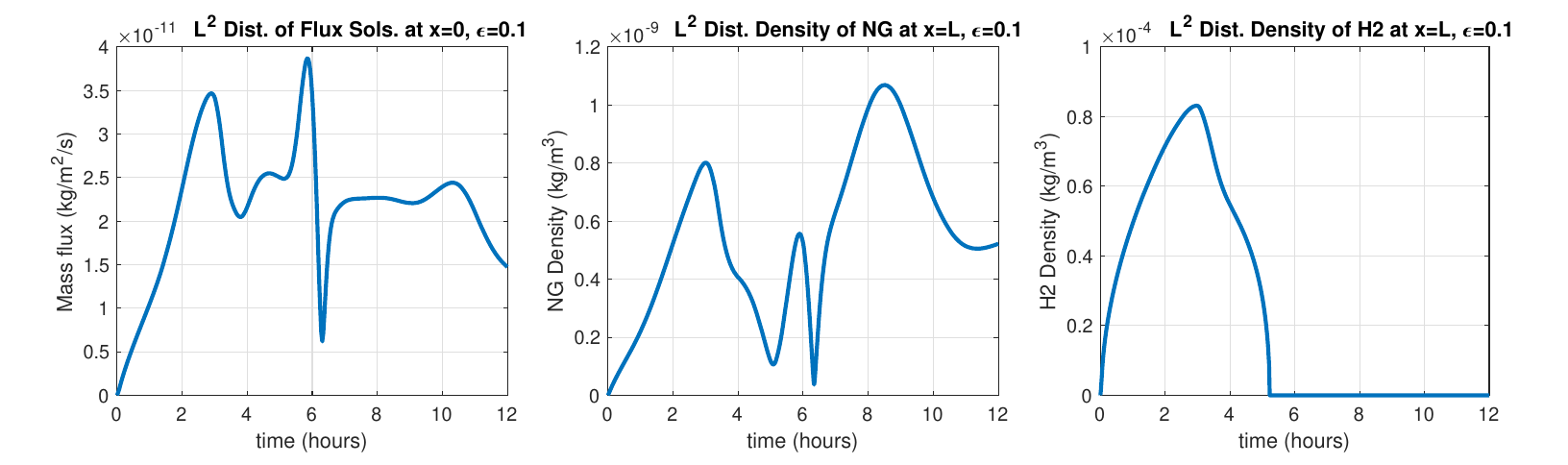}
    \vspace{-4ex}
	  \caption{Plots of $L^2$ differences in simulations with and without accounting for diffusion with  coefficients $\epsilon_1=\epsilon_2=0.1$.  Left: difference in input flux; Center: natural gas densities; Right: hydrogen densities.}
       \label{fig:one_pipe_mixing_mismatch_01}
\end{figure}

\begin{figure}[ht!]
    \centering
     \vspace{-1ex}
	  \includegraphics[width=1\linewidth]{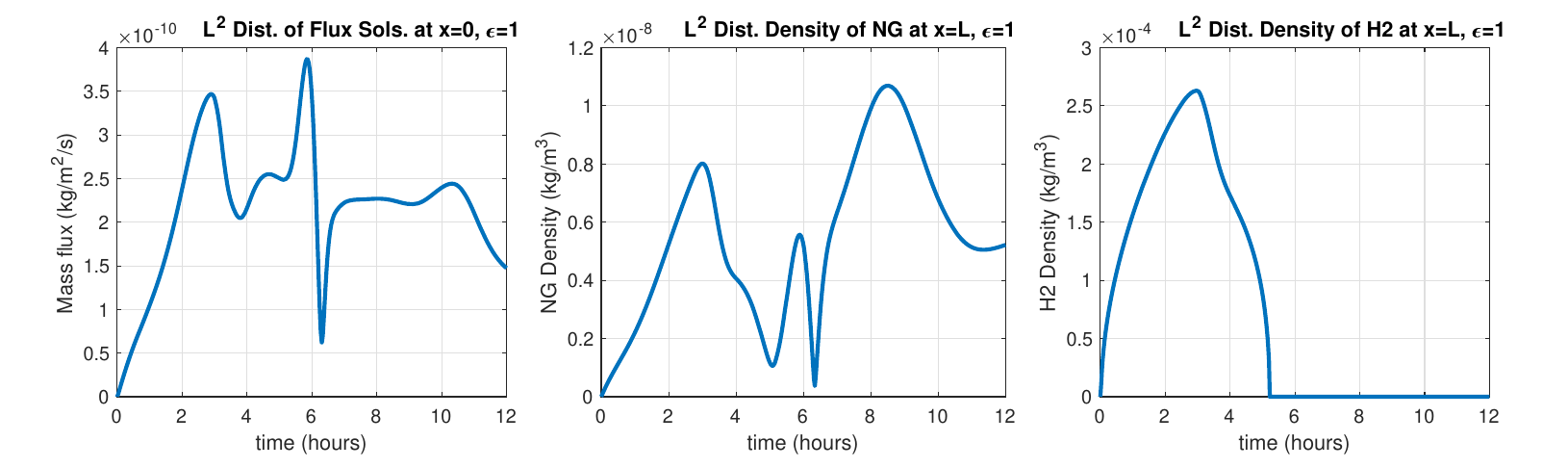}
    \vspace{-4ex}
	  \caption{Plots of $L^2$ differences in simulations with and without accounting for diffusion with  coefficients $\epsilon_1=\epsilon_2=1$.  Left: difference in input flux; Center: natural gas densities; Right: hydrogen densities.}
       \label{fig:one_pipe_mixing_mismatch_1}
\end{figure}

\noindent Observe that the simulation outputs shown in Figure \ref{fig:one_pipe_one_gas} exactly match the solution of the original IBVP \cite{herty2010new}. Next, we examine simulations of the problem with $c(t)$ as given in the definition above, where the results shown in Figure \ref{fig:one_pipe_mixing} are produced assuming that the diffusion coefficients are $\epsilon_1=\epsilon_2=0$. We then investigate the significance of diffusion effects in the dynamics by performing simulations with various values of the diffusion coefficients, and then plot and inspect the differences with respect to the baseline simulation shown in Figure \ref{fig:one_pipe_mixing}.  The $L^2$ differences in simulations with and without accounting for diffusion where values of $\epsilon_1=\epsilon_2=0.0001$, $\epsilon_1=\epsilon_2=0.1$, and $\epsilon_1=\epsilon_2=1$ are used are shown in Figures \ref{fig:one_pipe_mixing_mismatch_00001}, \ref{fig:one_pipe_mixing_mismatch_01}, and \ref{fig:one_pipe_mixing_mismatch_1}, respectively.

We observe by comparing the simulation discrepancies shown in Figures \ref{fig:one_pipe_mixing_mismatch_00001}, \ref{fig:one_pipe_mixing_mismatch_01}, and \ref{fig:one_pipe_mixing_mismatch_1} to the simulation solutions in Figure \ref{fig:one_pipe_mixing} that the differences between simulations that do and do not account for diffusion are many orders of magnitude smaller than the simulation outputs.  Because this holds for a wide range of values of diffusion coefficients, we suppose that we can omit the diffusion terms in the above problem formulation and set the right-hand side of the conservation of mass equations to zero.  To further justify this simplification, we direct the reader to a study on the empirical measurement of diffusion coefficient values for tracer gases in natural gas pipelines \cite{takeuchi2009flowmetering}, where the authors observe values that are within the range of our numerical experiments.

\subsection{Test Network Simulations} \label{sec:testnetwork_example}

To demonstrate our numerical method for general pipeline network topologies, we apply it to an IBVP that was examined in several previous studies \cite{gyrya2019explicit,bermudez2021modelling}.  We also cross-validate the method with solutions obtained in those studies, as well as with a recent lumped-element simulation method \cite{baker2023optimal}. To make our study as self-contained as possible, we define the network topology and parameters as well as the initial and boundary conditions for the benchmark IBVP.  We then conduct four comparison simulations.  First, we solve this IBVP for a single gas using our new staggered grid method for gas mixtures as well as the previous staggered grid simulation scheme for a single gas \cite{gyrya2019explicit}.  Second, we add a hydrogen injection and solve the modified IBVP using our new staggered grid method as well as a recently developed lumped element method \cite{baker2023optimal}.  Third, we compare two simulations of the modified IBVP using the staggered grid method with ideal and non-ideal gas equation of state modeling.  Finally, we conduct two simulations using our new method in which we simulate the modified IBVP both with and without the nodal monitoring policy.  All simulations are implemented in MATLAB on Dell G515(intel core i5 (8th Gen)).  We do not benchmark simulation time for this study, because our implementation is not optimized for computational efficiency.

\subsubsection{Test Network Baseline IBVP} \label{sec:testnetwork_def}

\begin{figure}[ht!]
    \centering
     \vspace{-1ex}
	  \includegraphics[width=1\linewidth]{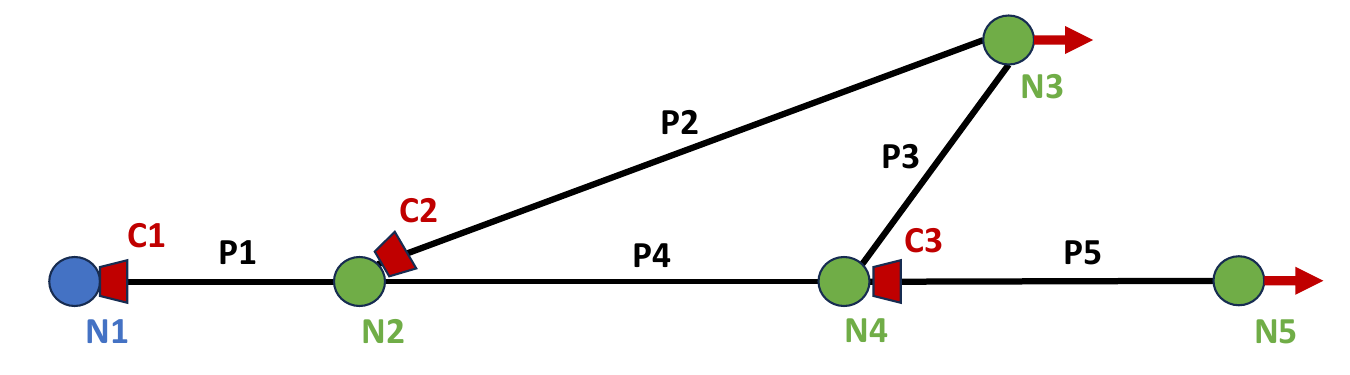}
    \vspace{-4ex}
    \caption{Test Network with five nodes,  five pipes, and three compressors.  The slack node N! is indicated in blue, and withdrawals at nodes N3 and N5 are indicated with red arrows.}
    \label{fig:network}
\end{figure}

The test network has five nodes, five pipes, and three compressors as illustrated in Figure \ref{fig:network}. The topology and parameters for pipes and compressors are given in Tables \ref{tab:pipe_physical_parameters} and \ref{tab:compressor_locations_in_network}, with enumeration as illustrated in the figure.
\begin{table}[h!]\textbf{}
    \centering
    \begin{tabular}{|| c | c | c | c | c | c ||} 
    \hline
    pipe ID & from node ID & to node ID & diameter (m) & length (km) & friction factor \\ [0.5ex] 
    \hline\hline
    P1 & N1 & N2 & 0.9144 & 20 & .01 \\ 
    \hline
    P2 & N2 & N3 & 0.9144 & 70 & .01 \\
    \hline
    P3 & N3 & N4 & 0.9144 & 10 & .01 \\
    \hline
    P4 & N2 & N4 & 0.6350 & 60 & .015 \\
    \hline
    P5 & N4 & N5 & 0.9144 & 80 & .01 \\
    \hline
    \end{tabular}
     \vspace{-1ex}
    \caption{Physical parameters for each pipe in the network.}
    \label{tab:pipe_physical_parameters}
\end{table}
%\vspace{-1ex}
\begin{table}[h!]
    \centering
    \begin{tabular}{|| c | c | c ||} 
    \hline
    comp ID & location node ID & to pipe ID \\ [0.5ex] 
    \hline\hline
    C1 & N1 & P1  \\ 
    \hline
    C2 & N2 & P2  \\
    \hline
    C3 & N4 & P5  \\
    \hline
    \end{tabular}
     \vspace{-1ex}
    \caption{Location of compressors in the pipe network.}
    \label{tab:compressor_locations_in_network}
\end{table}
The initial conditions for the IBVP are specified as steady-state flow based on the nodal boundary data given in Tables \ref{tab:network_initial_data_by_node} and \ref{tab:network_initial_data_by_pipe}.  The flow and endpoint pressures on each pipe can then be used to compute the gas densities on each edge according to the initial condition specification in the problem statement at the beginning of the section.
\begin{table}[ht!]
    \centering
    \begin{tabular}{|| c | c | c | c | c ||} 
    \hline
    item ID & item type & value type & variable & value \\ [0.5ex] 
    \hline\hline
    N1 & node & pressure (MPa) & $p_1(0)$ & 3.447378645 \\ 
    \hline
    N2 & node & flow withdrawal (kg/s) & $F_2^d(0)$ & 0 \\
    \hline
    N3 & node & flow withdrawal (kg/s) & $F_3^d(0)$ & 150 \\
    \hline
    N4 & node & flow withdrawal (kg/s) & $F_4^d(0)$ & 0 \\
    \hline
    N5 & node & flow withdrawal (kg/s) & $F_5^d(0)$ & 150 \\ 
    \hline
    C1 & comp & boost ratio & $\mu_1(0)$ & 1.5290113 \\
    \hline
    C2 & comp & boost ratio & $\mu_2(0)$ & 1.1128863 \\
    \hline
    C3 & comp & boost ratio & $\mu_3(0)$ & 1.2242249 \\ 
    \hline
    \end{tabular}
     \vspace{-1ex}
    \caption{Network initial data by node}
    \label{tab:network_initial_data_by_node}
\end{table}
\begin{table}[ht!]
    \centering
    \begin{tabular}{|| c | c | c | c ||} 
    \hline
    pipe ID & pressure in (MPa) & pressure out (Pa) & flow (kg/s) \\ [0.5ex] 
    \hline\hline
    P1 & 5.2710811 & 4.6112053 & 300.0 \\ 
    \hline
    P2 & 5.1317472 & 3.5400783 & 233.3 \\
    \hline
    P3 & 3.5400783 & 3.5043953 & 83.33 \\
    \hline
    P4 & 4.6112053 & 3.5043953 & 66.66 \\
    \hline
    P5 & 4.2901680 & 3.4473786 & 150.0 \\
    \hline
    \end{tabular}
     \vspace{-1ex}
    \caption{Network initial data by pipe.}
    \label{tab:network_initial_data_by_pipe}
\end{table}
The boundary conditions that define the evolution of the pressures and flows in the system are defined as transient nodal flow withdrawals $F_3^d(t)$ and $F_5^d(t)$ at nodes N3 and N5, respectively, and the compression ratios $\mu_1(t)$, $\mu_2(t)$, and $\mu_3(t)$ of the three compressors.  The mass withdrawal flows are given by
\begin{align}
    & F_3^d(t ~ \mbox{mod} ~ T) = \phi_3^d(0) \cdot \left( 1 - \frac{1}{10} \cdot \left( 1 - \cos{\frac{2\pi}{T} t} \right) \right),\\
    & F_5^d(t ~ \mbox{mod} ~ T) = \phi_5^d(0) \cdot \begin{cases}
        1, \quad 0<t\leq 12000\\
        \frac{1}{3} + \frac{5t}{90000}, \quad 12000<t\leq 15600\\
        1.2, \quad 15600<t\leq 48000,\\
        \frac{193}{50} - \frac{5t}{90000}, \quad 48000<t\leq 51600\\
        1, \quad 51600<t\leq 86400\\
    \end{cases}
\end{align}
\noindent where $T = 86400$ sec, and $F_3^d(0)$ and $F_3^d(0)$ are given in Table \ref{tab:network_initial_data_by_node} . Observe that $F_5^d(t)$ defines linear interpolation of the points $(0, F_5^d(0))$, $(12000, F_5^d(0))$, $(15600, 1.2 \cdot F_5^d(0))$, $(48000, 1.2 \cdot F_5^d(0))$, $(51600, F_5^d(0))$, and $(86400, c_2(0))$.  The time-varying compressor ratios are given by

\begin{align*}
    & \mu_1(t ~ \mbox{mod} ~ T) = \mu_1(0) \cdot \left( 1 - \frac{1}{10} \cdot \left( 1 - \cos{\frac{2\pi}{T} t} \right) \right),\\
    & \mu_2(t ~ \mbox{mod} ~ T) = \mu_2(0) \cdot \begin{cases}
        1, \quad 0<t\leq 21600\\
        -1.4 + \frac{t}{9000}, \quad 21600<t\leq 25200\\
        1.4, \quad 25200<t\leq 64800,\\
        8.6 - \frac{t}{9000}, \quad 64800<t\leq 68400\\
        1, \quad 68400<t\leq 86400\\
    \end{cases} \\
    & \mu_3(t ~ \mbox{mod} ~ T) = \mu_3(0) \cdot \left( 1 + \frac{1}{4} \cdot \left( 1 - \cos{\frac{6\pi}{T} t} \right) \right),
\end{align*}

\noindent where $T = 86400$ sec, and $\mu_1(0)$, $\mu_2(0)$, and $\mu_3(0)$ are given in Table \ref{tab:network_initial_data_by_node} . Observe that $\mu_2(t)$ defines linear interpolation of the points $(0, \mu_2(0))$, $(21600, \mu_2(0))$, $(25200, 1.4 \cdot \mu_2(0))$, $(64800, 1.4 \cdot \mu_2(0))$, $(68400, \mu_2(0))$, and $(86400, \mu_2(0))$.  The initial values of the boundary conditions are equal to the initial nodal values, and the resulting IBVP is well-posed.  Note that the supply at node N1 is homogeneous natural gas, and the above modeling does not include any hydrogen.  The injection of hydrogen will be specified subsequently for the simulations in Section \ref{sec:testnetwork_mixture}.
\begin{figure}[h!]
    \centering
     \vspace{-1ex}
    \includegraphics[width=1\linewidth]{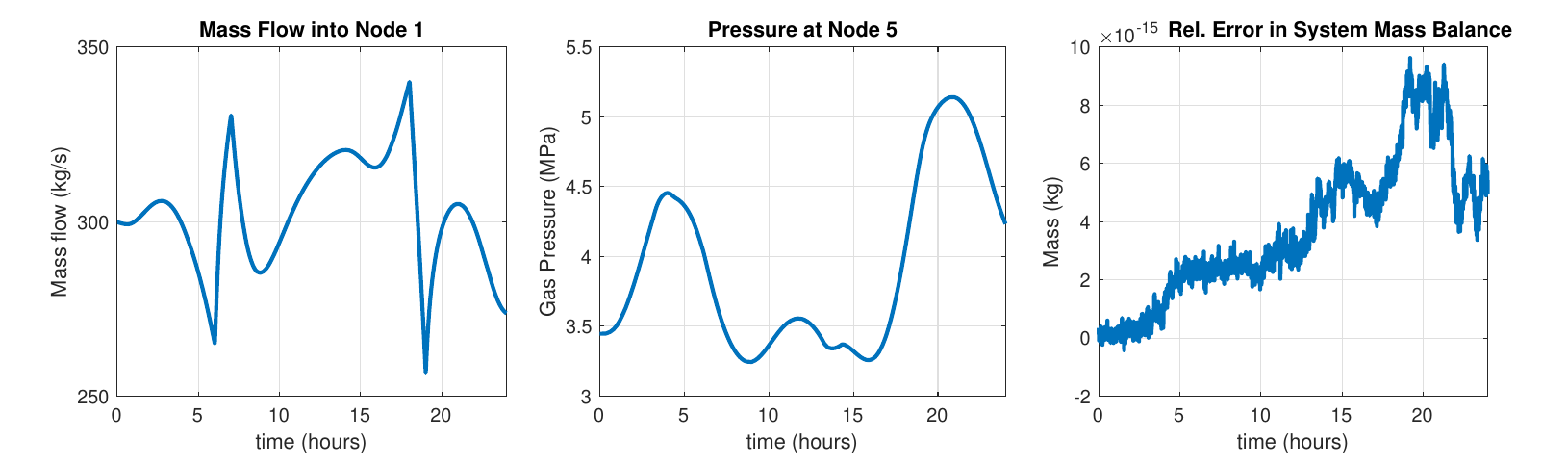}
     \vspace{-4ex}
    \caption{Simulation of the network IBVP defined above using the ideal gas equation of state. Left: Inflow to node N1; Center: pressure at node N5; Right: the discrepancy between the change in mass of the gas in the pipe and the mass balance of gas entering and leaving the system through the nodes. The error in mass conservation is on the order of machine precision.}
    \label{fig:one gas}
\end{figure}

\subsubsection{Validation for a Single Gas}  \label{sec:testnetwork_onegas}

In order to verify our new staggered grid finite difference scheme for heterogeneous gas mixing in flow through a pipeline network, we first simulate the IBVP defined above in Section \ref{sec:testnetwork_def} for homogeneous natural gas, and compare the results with the previously developed staggered grid scheme that accommodates only a single mass conservation equation \cite{gyrya2019explicit}. Figure \ref{fig:one gas} shows a comparison of the solutions obtained by the two schemes using a time discretization of $\Delta t = .02$. We also examine the conservation of mass within the pipe by integrating the density along the length of the pipe and subtracting the difference between flows entering and leaving the pipe.  This conservation of mass is verified for both methods.

\begin{figure}[ht!]
    \centering
     \vspace{-2ex}
    \includegraphics[width=1\linewidth]{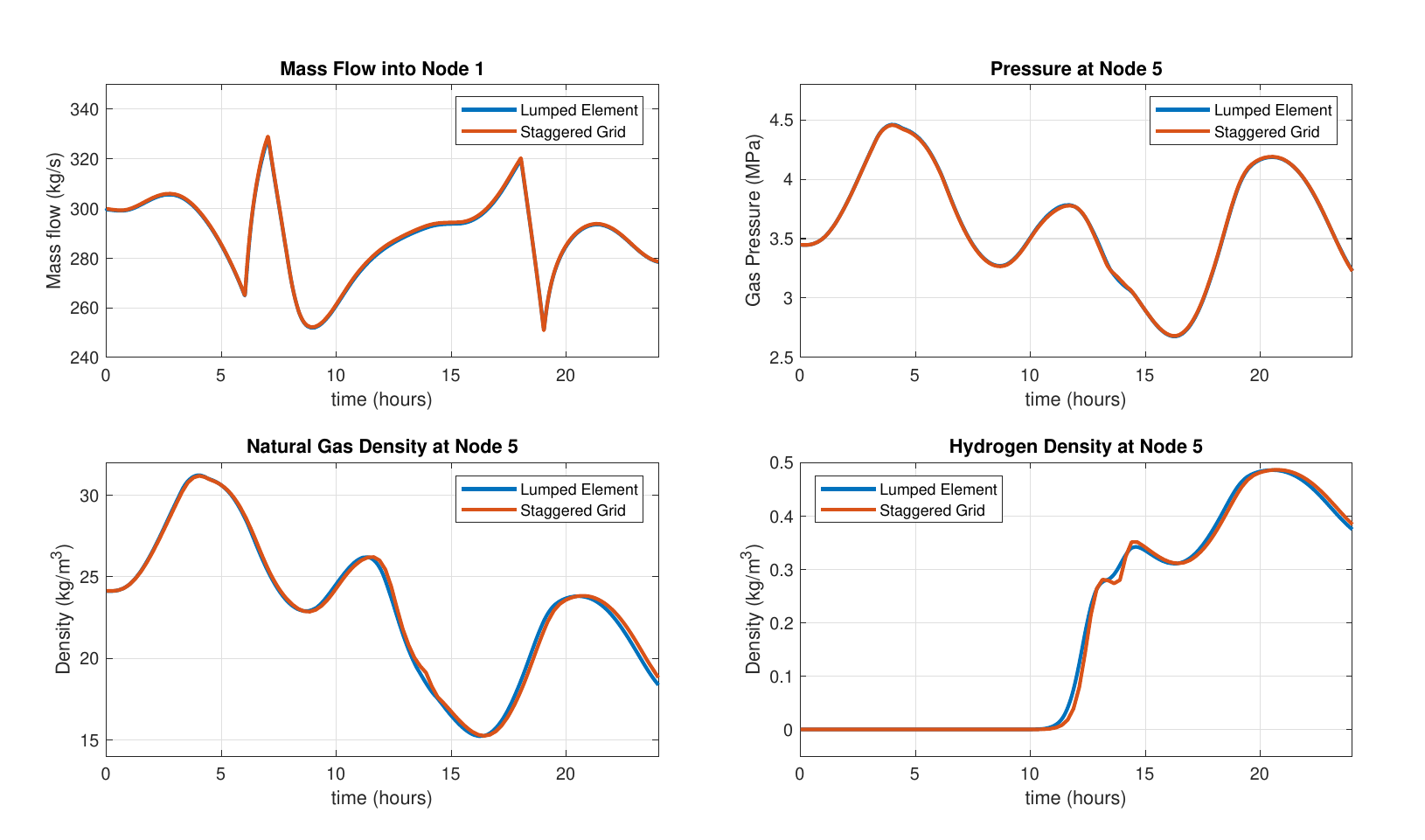} 
     \vspace{-4ex}
    \caption{Comparison of staggered grid discretization using $\Delta t=0.02$ seconds and lumped element method for the IBVP defined in Section \ref{sec:onepipe_example} with source hydrogen concentration given in equation \eqref{eq:hydrogen_condition_example}.  The ideal gas model is used. Top Left: mass inflow to node N1; Top Right: pressure at node N5; Bottom Left: density of natural gas at node N5; Bottom Right: density of hydrogen at node N5.}
    \vspace{-1ex}
    \label{fig:comparison_with_lumped_elements}
\end{figure}
\subsubsection{Validation for a Gas Mixture}  \label{sec:testnetwork_mixture}

A compelling means of verifying that a computational method is correct is by conducting a cross-comparison of a different computational method for the same IBVP.  In Figure \ref{fig:comparison_with_lumped_elements}, we show the results of two simulations of the IBVP defined in Section \ref{sec:onepipe_example}, where one is conducted using our staggered grid scheme and another set of results is obtained using the lumped element model described in a recent study \cite{baker2023optimal}.  The additional hydrogen injection occurs at node N1, and when $T=86400$ seconds takes the form
\begin{equation} \label{eq:hydrogen_condition_example}
    c_1^{\alpha}(t)=0.01\cdot (1+\tanh{(0.0005\cdot(t-T/3)))}.
\end{equation}
The staggered grid simulation was implemented directly using MATLAB on Dell G515 with Intel core i5 (8th Gen) processor using a time discretization $\Delta t = .02$.  The lumped element simulation was done on a MacBook Air 8-core CPU with 8GB of unified memory, and is implemented in MATLAB using the function using the function \verb"ode15s".

\subsubsection{Comparison of Ideal and Non-Ideal Gas Modeling}  \label{sec:testnetwork_eos}

\begin{figure}[ht!]
    \centering
     \vspace{-2ex}
    \includegraphics[width=1\linewidth]{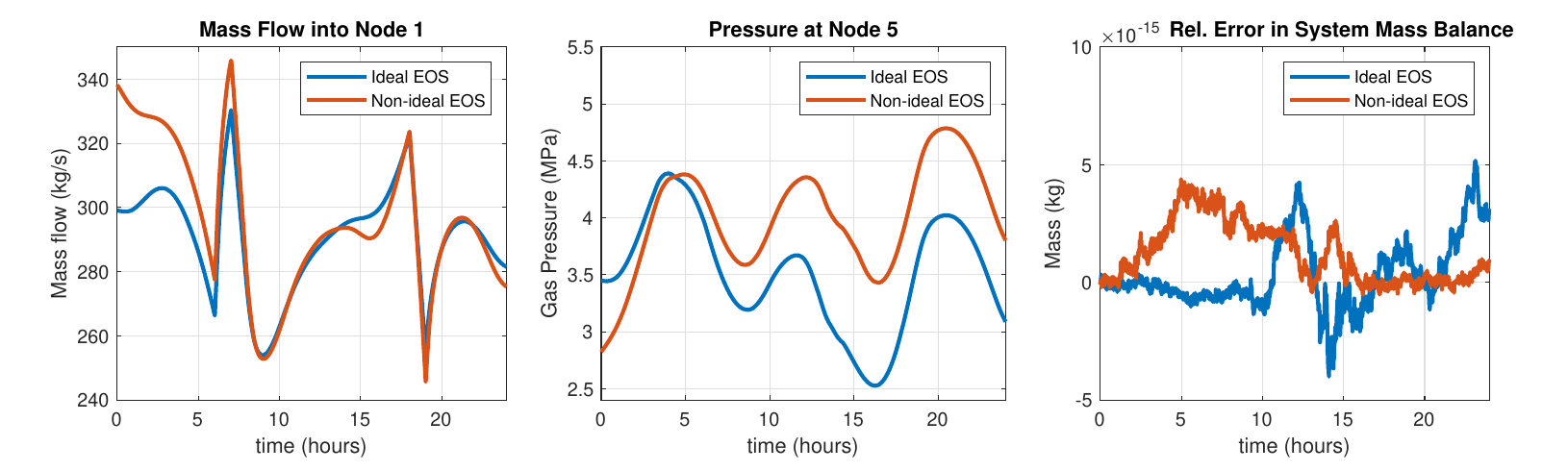} 
     \vspace{-4ex}
    \caption{Comparison of ideal and non-ideal gas modeling for the staggered grid discretization simulation of mixing of gases in the pipeline network using time step $\Delta t=0.1$. Left: inflow through node N1; Center: pressure at node N5; Right:  discrepancy between the mass of the gas in the pipe and the mass of the gas leaving/entering the system through the nodes.}
    \label{fig:ideal_nonideal}
\end{figure}

We now contrast two simulations of heterogeneous gas flow defined by the IBVP in Section \ref{sec:onepipe_example} to highlight the significance of appropriate equation of state modeling.  Hydrogen is injected at node N1 according to the time-varying profile in equation \eqref{eq:hydrogen_condition_example}.   We use the equation of state modeling presented in Section \ref{sec:EOS} to compare a simulation using ideal gas modeling and nonideal gas modeling where the compressibility factors are approximated by linear functions of pressure in the form
\begin{align*}
	%& Z_{NG}(p) \approx 1 - 0.000000025 p,\\
	Z_{NG}(p) &\approx 1 - (0.25 \times 10^{-7}) \cdot p,\\
	%& Z_{H}(p) \approx  1 + 0.0000000059 p.
	Z_{H}(p) &\approx  1 + (0.59 \times 10^{-8}) \cdot p.
\end{align*}
The comparison of simulation results is shown in Figure \ref{fig:ideal_nonideal}, in addition to the mass balance verification for both simulations in the same sense as shown in Section \ref{sec:testnetwork_onegas}.  By inspection, it is evident that appropriate gas equation of state modeling makes a significant difference in simulation results, and is critical to ensure that the monitoring policies proposed in Section \ref{sec:monitoring} and demonstrated below in Section \ref{sec:testnetwork_monitoring} trigger corrective actions at the correct pressure and/or hydrogen fraction levels.

\subsubsection{Demonstration of the Nodal Monitoring Policy}  \label{sec:testnetwork_monitoring}

Finally, we demonstrate the implementation of the nodal monitoring policy for the maximum mass fraction described in Section \ref{sec:input_flow_monitoring_concentration}.  
We consider an ideal gas model for simplicity, with hydrogen mass fraction for the supply at node N1 specified by equation \eqref{eq:hydrogen_condition_example}.  Furthermore, we add a constant injection flow of 2 \kg/s hydrogen at node N4 throughout the duration of the 24 hour simulation period.  

\begin{figure}[h!]
    \centering
     \vspace{-2ex}
    \includegraphics[width=1\linewidth]{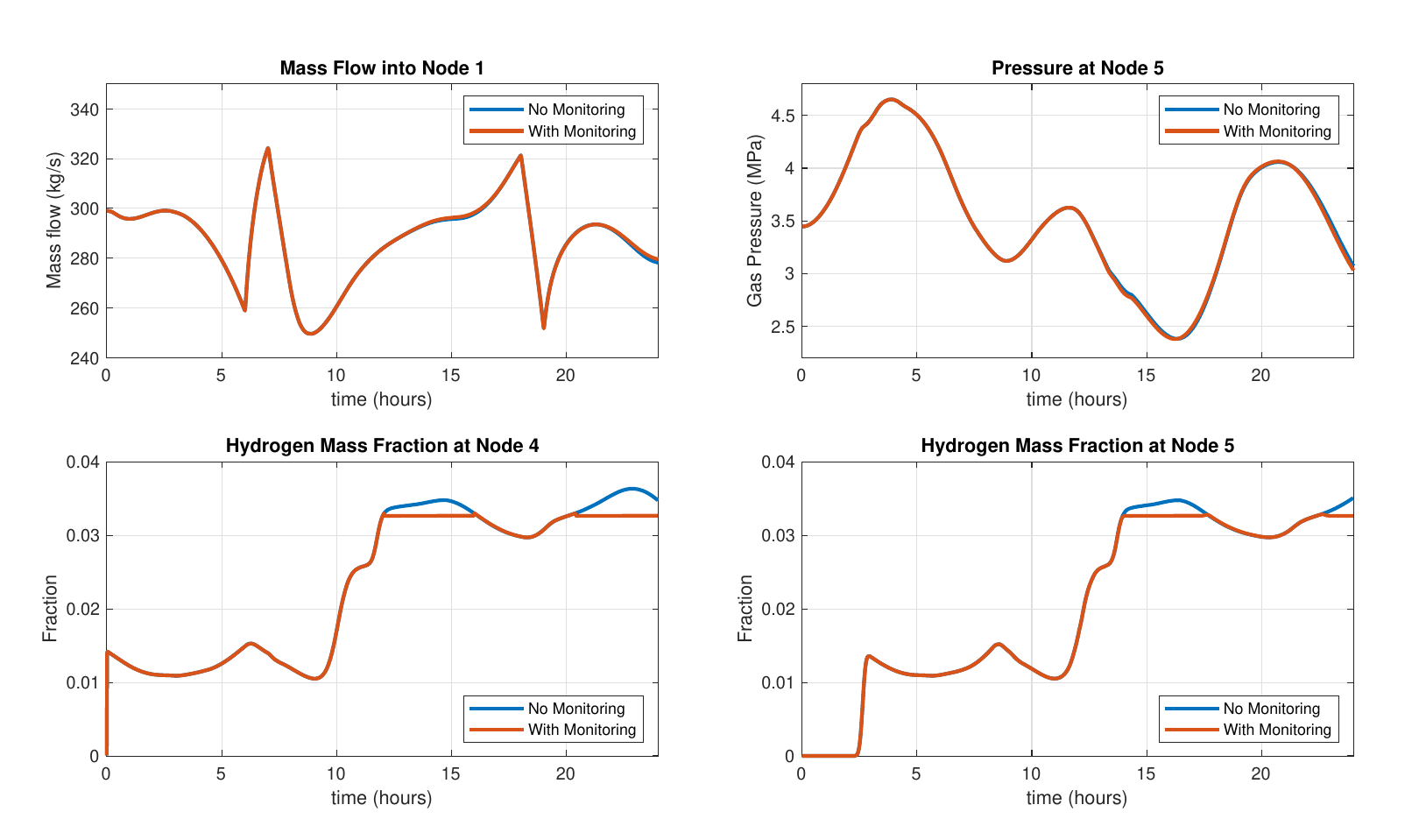} 
     \vspace{-4ex}
    \caption{Demonstration of nodal monitoring. Ideal gas modeling and $\Delta t=0.1$ are used. Top Left: inflow through node 1; Top Right: pressure at node 5; Bottom Left: hydrogen mass fraction at node 4; Bottom Right: hydrogen mass fraction at node 5. }
    \label{fig:monitoring}
\end{figure}

The monitoring policy is applied using an upper bound value for the hydrogen mass fraction of 3.3\%, and the results are shown in Figure \ref{fig:monitoring}.  By inspecting the figure, one sees that the implementation of the nodal monitoring policy to control the hydrogen injection at node N4 prevents the hydrogen mass fraction from exceeding the prescribed limit.  Without the policy, we observe clear violations of the constraint.
Both simulations were made using MATLAB on Dell G515 (Intel Core i5 (8th Gen)) using a time discretization $\Delta t = .1$ seconds. 

%\newpage

\section{Conclusions}  \label{sec:conclusions}

We have presented the first explicit second order staggered finite difference discretization scheme for simulating the transport of highly heterogeneous gas mixtures through pipeline networks.  The method can be used to accurately simulate the blending of hydrogen into natural gas pipelines to reduce end use carbon emissions while using existing pipeline systems throughout their planned lifetimes.  Our computational method accommodates an arbitrary number of constituent gases with widely varying physical properties that may be injected into a network with significant spatiotemporal variation, which makes the gas flow physics highly location- and time- dependent.  Notably, the accommodation of non-ideal equations of state enables modeling the propagation of pressure dynamics based on locally variable wave speeds that depend on mixture composition and density.  We derive compatibility relationships for network edge domain boundary values that are significantly more complex than in the case of a homogeneous gas.   The key innovation of our method is the fully explicit computation of all simulation steps, without the need to resort to implicit calculations.  By precisely accounting for the local composition and nodal mixing, our scheme can be used as a benchmark for validating coarser models for optimization and sensitivity studies.  

\newpage

\section{APPENDIX} \label{sec:appendix}

We provide several intermediate derivations to assist the reader to interpret the results in Section \ref{sec:Nodal_BC} and Section \ref{sec:monitoring} on nodal boundary conditions and nodal monitoring policies.  We also provide details on the values used for the coefficients $a_\alpha$ in the linear approximation of gas compressibility factors in equation \eqref{eq:EOS_linear_Z}. 

\subsection{Nodal boundary conditions} \label{sec:append_nodal_bc}

Recall equations \eqref{eq:disc_flow_balance_at_node} and \eqref{eq:conservation_of_momentum_discrete_alt} for flow balance at a node $q\in\mathcal V$ and conservation of momentum at the start ($j=0$) of a pipe $k\in\mathcal{E}$ leaving node $q$, respectively:
\begin{equation}\label{eq:disc_flow_balance_at_node_appendix}
    \sum_{k\in\partial^+q} S_k \phi_{m,j=N_k}^{\a,k}-\sum_{k\in\partial^-q} S_k \phi_{m,j=0}^{\a,k} = F_m^{\a,q,d} - F_m^{\a,q,s}, \quad \forall \,\alpha\in\{1,\ldots,n\}, \, \forall \, q\in\mathcal{V},
\end{equation}
\begin{align}\label{eq:conservation_of_momentum_discrete_alt_appendix}
	\frac{1}{\dt}\left(\phi_{m,j=0}^{k} - \phi_{m-1,j=0}^{k}\right) 
	+ 
	\frac{1}{\dx}
	\left(p^{k}_{n,i=1/2} - \mu_n^k p_n^{q} \right) 
	= 
	-
	f(\phi_{n,j=0}^{k},d_{n,j=0}^{k}).
\end{align}
Recall also that we may write equation \eqref{eq:disc_flow_balance_at_node_appendix} in a simplified form, where we supposed that the flow directions are all oriented out of node $q$ and signed appropriately, as
\begin{equation}\label{eq:disc_flow_balance_at_node_simple_appendix}
    \sum_{k\in\partial q} S_k \phi_{m,j=0}^{\a,k} + F_m^{\a,q,d} - F_m^{\a,q,s} = 0, \quad \forall \,\alpha\in\{1,\ldots,n\}, \, \forall \, q\in\mathcal{V}.
\end{equation}
We suppose that the pressure is as defined in equation \eqref{eq:explicit_pressure} as
\begin{align}\label{eq:pressure_def_appendix}
	p^{k}_{n,i} = p(\mathbf{d}^{k}_{n,i}), 
\end{align}
and the friction term is
\begin{align}\label{eq:friction_discrete_alt_2_appendix}
	f(\phi_{n,j=0}^{k},d_{n,j=0}^{k})
	=
	\frac{\lambda_k}{2D_k} \cdot
	\frac{\phi_{n,j=0}^{k}|\phi_{n,j=0}^{k}|}{d_{n,j=0}^{k}}.
\end{align}
Recall also that $c^\alpha=d^\alpha/d$ is the mass fraction of gas $\alpha$ where $d=\sum_\alpha d^\alpha$ is the total density, so that mass fraction at a node $q$ and on a pipe $k$ is defined at time $m$ as
\begin{equation} \label{eq:concentration_as_density_fraction_appendix}
    c_m^{\a,q}=\frac{d_m^{\alpha,q}}{\sum_{\beta} d_m^{\beta,q}}, \,\, \forall q\in\mathcal{V}, \quad c_{m,i}^{\a,k}=\frac{d_{m,i}^{\alpha,k}}{\sum_{\beta} d_{m,i}^{\beta,k}}, \,\, \forall k\in\mathcal{E}.
\end{equation}
We then express the gas component mass flux at the start $j=0$ of a pipe $k$ at time $m$ as 
\begin{equation} \label{eq:component_flux_appendix}
    \phi_{m,j=0}^{\a,k}=\phi^{k}_{m,j=0}\cdot c_{m-1/2, i = 1/2}^{\a, k} , \,\, \forall k\in\mathcal{E},
\end{equation}
where mass fraction is computed a half-step adjacent on the staggered grid.
Using the above definitions, we can explicitly state the nodal pressure $p^q_n$ at node $q$ and time $n$, and the boundary flow $\phi^{k}_{m,j=0}$ for pipe $k$ at time $m$ and densities $d_m^{\alpha,q}$ at node $q$ as functions of nodal pressure $p^q_n$. 

When using ideal gas equation of state modeling for each separate pipe, we solve equation \eqref{eq:conservation_of_momentum_discrete_alt_appendix} for the boundary flow $\phi^{k}_{m,j=0}$, where the only unknown quantity is $p^q_n$, so that we write the dependence as $\phi^{k}_{m,j=0} = \phi^{k}_{m,j=0}\left(p^q_n\right)$: 	
\begin{equation}\label{eq:conservation_of_momentum_discrete_flux_pipe_ideal_appendix}
	\phi^{k}_{m,j=0}\left(p^q_n\right)
	\!= \! \left(\!\phi_{m-1,j=0}^{k}
	\!-\!
	\frac{\dt \lambda_k}{2D_k}
	\frac{\phi_{m-1,j=0}^{k}|\phi_{m-1,j=0}^{k}|}{\sum_\a d^{\a,k}_{m - 1/2,i = 1/2}} \!-\! \frac{\dt}{\dx} \sum_\a d^{\a,k}_{m - 1/2,i = 1/2} R_\a T \!\right) \!+\! \frac{\dt}{\dx} \mu_n^k p^q_n, \!
\end{equation}
where $\mu_n^k\equiv 1$ if pipe $k$ has no compressor at the start.  Similarly, when using non-ideal gas equation of state modeling for each separate pipe, we solve equation \eqref{eq:conservation_of_momentum_discrete_alt_appendix} for the boundary flow $\phi^{k}_{m,j=0}$ as a function of the unknown quantity $p^q_n$ to obtain:	
\begin{align}\label{eq:conservation_of_momentum_discrete_flux_pipe_nonideal_appendix}
	\phi^{k}_{m,j=0}\left(p^q_n\right)
	\!=\! \left(\!\phi_{m-1,j=0}^{k}
	\!-\!
	\frac{\dt \lambda_k}{2D_k}
	\frac{\phi_{m-1,j=0}^{k}|\phi_{m-1,j=0}^{k}|}{\sum_\a d^{\a,k}_{m - 1/2,i = 1/2}} \!-\! \frac{\dt}{\dx} \frac{\sum_\a d^{\a,k}_{m - 1/2,i = 1/2} R_\a T}{1 \!-\! \sum_\a d^{\a,k}_{m - 1/2,i = 1/2} R_\a T a_\a} \!\right) \!\!+\! \frac{\dt}{\dx} \mu_n^k p^q_n,
\end{align}
where $\mu_n^k\equiv 1$ if pipe $k$ has no compressor at the start. 
Next, let us examine the equation system for nodal flow balance conditions \eqref{eq:disc_flow_balance_at_node_appendix} for all gas constituents $\alpha=\{1,\ldots,n\}$.  Given any gas withdrawals from and injections into the network at node $q$ at time point $m$ denoted by $F_m^{\alpha,q,d}$ and $F_m^{\alpha,q,s}$, respectively, for each gas $\alpha$, we have
\begin{equation}\label{eq:flow_balance_vector_1_appendix}
		\left(\sum_{k\in\partial_+r} S_k \phi^{k}_{m,j=N_k}\left(p^q_n\right) c_{m-1/2, i = N_k-1/2}^{\a, k}  \right) + F_m^{\alpha,q,s}
	=  F_m^{\alpha,q,d} + \sum_{k\in\partial_-r}  S_k \phi^{k}_{m,j=0}\left( p^q_n\right)  \cdot
		\frac{d_m^{\a,q}}{\sum_\beta d_m^{\beta,q}}, \,\, \forall \a.
\end{equation}
Alternatively, equation \eqref{eq:flow_balance_vector_1_appendix} can be written entirely in terms of mass fractions as
\begin{equation}\label{eq:flow_balance_vector_2_appendix}
		\left(\sum_{k\in\partial_+r} S_k \phi^{k}_{m,j=N_k}\left(p^q_n\right) c_{m-1/2, i = N_k-1/2}^{\a, k}  \right) + F_m^{\alpha,q,s}
	= \left(F_m^{q,d} + \sum_{k\in\partial_-r}  S_k \phi^{k}_{m,j=0}\left( p^q_n\right) \right) \cdot
		c_m^{\a,q},  \,\, \forall \a.
\end{equation}
We can then solve equation \eqref{eq:flow_balance_vector_2_appendix} to obtain expressions for individual nodal mass fractions and densities after mixing:		
\begin{align}%\label{eq:explicit_nodal_density2_appendix}
		c_m^{\alpha,q} & = \frac{\sum_{k\in\partial_+r} S_k \phi^{k}_{m,j=N_k}\left(p_n^q\right)c_{m-1/2, i = 1/2}^{\a,k} + F_m^{\alpha,q,s} }{F_m^{q,d} + \sum_{k\in\partial_-r}  S_k \phi^{k}_{m,j=0}\left( p_n^q\right) }, \label{eq:explicit_nodal_density2_appendix_c} \\
        d_m^{\alpha,q} & = \frac{\sum_{k\in\partial_+r} S_k \phi^{k}_{m,j=N_k}\left(p_n^q\right)c_{m-1/2, i = 1/2}^{\a,k} + F_m^{\alpha,q,s} }{F_m^{q,d} + \sum_{k\in\partial_-r}  S_k \phi^{k}_{m,j=0}\left( p_n^q\right) } \cdot \sum_\beta d_m^{\beta,q}.  \label{eq:explicit_nodal_density2_appendix_d}
\end{align}
Recall now that substituting \eqref{eq:conservation_of_momentum_discrete_alt_appendix} into the summation of \eqref{eq:disc_flow_balance_at_node_simple_appendix} over $\alpha$ and reformulating yields an explicit expression (eq. \ref{eq:flow_balance_scalar_coef1}) for the nodal pressure $p_n^{q}$ at a node $q$ at time $t_n$, given by 
\begin{align}\label{eq:flow_balance_scalar_coef1_appendix}
	& p_n^{q} \!=\! \frac{\dx}{\dt}\!\cdot\! \frac{ F^{q,s}_m-F^{q,d}_m - \sum_{k\in\partial q} S_k  \left(\phi_{m-1,j=0}^{k} - \frac{\dt \lambda_k}{2D_k} \cdot \frac{\phi_{m-1,j=0}^{k}|\phi_{m-1,j=0}^{k}|}{\sum_\a d_{m - 1/2,i = 1/2}^{\a,k}} - \frac{\dt}{\dx} p_{n,i=1/2}^{k} \right)}{\sum_{k\in\partial q} S_k\mu_n^k},
\end{align}	
where $\mu_n^k\equiv 1$ for all pipes $k$ with no compressors.  
	
When using the ideal gas approximation, we suppose that compressibility for each gas $\alpha$ is $Z_\alpha\equiv 1$ and apply additive partial pressures to express $p_{n,i=1/2}^{k}$ in equation \eqref{eq:flow_balance_scalar_coef1_appendix} in terms of partial densities $d^{\alpha}$ following equation \eqref{eq:EOS_linear_compressibility}, resulting in
\begin{equation} \label{eq:nodal_pressure_explicit_ideal_pres_in_dens_appendix}
p_{n,i=1/2}^{k}=\sum_\a d^{\a,k}_{m - 1/2,i = 1/2} R_\a T.
\end{equation}
When using non-ideal gas modeling with linear approximations of individual gas compressibility factors, we apply equation \eqref{eq:explicit_pressure} to obtain the expression for $p_{n,i=1/2}^{k}$ to use in equation \eqref{eq:flow_balance_scalar_coef1_appendix}, resulting in
\begin{equation} \label{eq:nodal_pressure_explicit_nonideal_pres_in_dens_appendix}
p_{n,i=1/2}^{k}=\frac{\sum_\a d^{\a,k}_{m - 1/2,i = 1/2} R_\a T}{1 - \sum_\a d^{\a,k}_{m - 1/2,i = 1/2} R_\a T a_\a}.
\end{equation}
Next, in order to compute component densities at the pipe boundaries, we require an expression for nodal densities in terms of nodal pressure.  Flux unwinding requires knowledge of nodal densities to determine boundary densities.  Inspecting equation \eqref{eq:explicit_nodal_density2_appendix_d}, we see that there is a single coefficient $g(p_n^q)$ depending on nodal pressure $p_n^q$ that relates the partial densities $d^{\alpha}_m$ to the left hand side of the equation:
\begin{align}
 d_m^{\alpha,q} & = \left(\sum_{k\in\partial_+r} S_k \phi^{k}_{m,j=N_k}\left(p_n^q\right)c_{m-1/2, i = 1/2}^{\a,k} + F_m^{\alpha,q,s} \right) \cdot g(p_n^q),  \quad \forall \alpha, \label{eq:explicit_nodal_density_appendix_d}\\
   g(p_n^q) & = \left(F_m^{q,d} + \sum_{k\in\partial_-r}  S_k \phi^{k}_{m,j=0}\left( p_n^q\right) \right)^{-1} \cdot \sum_\beta d_m^{\beta,q}. \label{eq:explicit_nodal_density_appendix_g}
\end{align}
The above expression for nodal densities can then be used together with the pressure-density relation \eqref{eq:explicit_pressure} in order to obtain a fully explicit form of the coefficient $g(p_n^q)$.

\subsubsection{Pipe Boundary Updates with Ideal Gas Modeling}

In the case of ideal gas modeling, we substitute equation \eqref{eq:explicit_nodal_density_appendix_d} into the formula \eqref{eq:EOS_linear_compressibility} for pressure-density dependence for an ideal gas to obtain		
\begin{align}\label{eq:nodal_pressure_ideal_gas_appendix}
	p_n^q
	= \sum_\a \left(\sum_{k\in\partial_+r} S_k \phi^{k}_{m,j=N_k}\left(p_n^q\right)c_{m-1/2, i = 1/2}^{\a,k} + F_m^{\alpha,q,s} \right) \!\cdot\! g(p_n^q) \!\cdot\! R_\a T .
\end{align}
We can rearrange equation \eqref{eq:nodal_pressure_ideal_gas_appendix} to obtain an explicit expression for the coefficient $g(p_n^q)$ as a function of the nodal pressure $p_n^q$:
\begin{align}\label{eq:g_coeff_as_fun_of_pres_ideal_appendix}
	g(p_n^q) = \frac{p_n^q}{\sum_\a \left(\sum_{k\in\partial_+r} S_k \phi^{k}_{m,j=N_k}\left(p_n^q\right)c_{m-1/2, i = 1/2}^{\a,k} + F_m^{\alpha,q,s} \right) R_\a T}.
\end{align}
Then substituting equation \eqref{eq:g_coeff_as_fun_of_pres_ideal_appendix} into equation \eqref{eq:explicit_nodal_density_appendix_d} leads to an expression for nodal component densities after mixing at node $q$ with linear dependence on nodal pressure $p_n^q$ and explicit dependence on boundary flows and densities.  We write this dependence as $d_n^{\a,q} = d_n^{\a,q}(p_n^q)$, with form
\begin{align}\label{eq:explicit_nodal_density_ideal_appendix}
		d_n^{\a,q}(p_n^q) = p_n^q \cdot \frac{\sum_{k\in\partial_+r} S_k \phi^{k}_{m,j=N_k}\left(p_n^q\right)c_{m-1/2, i = 1/2}^{\a,k} + F_m^{\alpha,q,s}}{\sum_\a \left(\sum_{k\in\partial_+r} S_k \phi^{k}_{m,j=N_k}\left(p_n^q\right)c_{m-1/2, i = 1/2}^{\a,k} + F_m^{\alpha,q,s} \right) R_\a T}.  
\end{align}
We can then formulate the time update of boundary flows for each pipe $k\in\mathcal{E}$ based on the dependence on nodal pressure  $p^{r,n}$ as	
\begin{equation}\label{eq:boundary_flow_as_f_of_p_ideal_appendix}
	\phi_{m,j=0}^{k}(p_n^q)
	\!=\! \left(\phi_{m-1,j=0}^{k} 	\!-\!
	\frac{\dt \lambda_k}{2D_k}
	\frac{\phi_{m-1,j=0}^{k}|\phi_{m-1,j=0}^{k}|}{\sum_\a d_{m - 1/2,i = 1/2}^{\a,k}} \!-\! \frac{\dt}{\dx} \sum_\a d_{m - 1/2,i = 1/2}^{\a,k} R_\a T \right) \!+\! \frac{\dt}{\dx} p_n^q.
\end{equation}
The time update for outgoing boundary densities for each gas component $\alpha$ for every pipe $k\in\mathcal{E}$ is obtained as
\begin{align}\label{eq:conservation_of_mass_discrete_ideal_appendix}
	d_{m + 1/2,i = 1/2}^{\a,k} = d_{m - 1/2,i = 1/2}^{\a,k}  
	- 
	\frac{\dt}{\dx}
	\left( \phi^{k}_{m,j=1} \frac{ d^{\a,k}_{m - 1/2,i = 1/2}}{\sum_\beta d^{\beta,k}_{m - 1/2,i = 1/2}} - \phi^{k}_{m,j=0}(p_n^q) \frac{d_n^{\a,q}(p_n^q)}{\sum_\beta d_n^{\beta,q}(p_n^q)} \right).
\end{align}
Critically, both updating expressions \eqref{eq:boundary_flow_as_f_of_p_ideal_appendix} and \eqref{eq:conservation_of_mass_discrete_ideal_appendix} have explicit dependence on data for previous time points on the staggered grid that are already known.

\subsubsection{Pipe Boundary Updates with Non-Ideal Gas Modeling}

In the case of non-ideal gas modeling, we substitute equation \eqref{eq:explicit_nodal_density_appendix_d} into the formula \eqref{eq:explicit_pressure} for pressure-density dependence for a non-ideal gas to obtain	
\begin{align}\label{eq:nodal_pressure_nonideal_gas_appendix}
	p_n^q
	=  \frac{\sum_\a \left(\sum_{k\in\partial_+r} S_k \phi^{k}_{m,j=N_k}\left(p_n^q\right)c_{m-1/2, i = 1/2}^{\a,k} + F_m^{\alpha,q,s} \right) \!\cdot\! g(p_n^q) \!\cdot\! R_\a T}{1 - \sum_\a \left(\sum_{k\in\partial_+r} S_k \phi^{k}_{m,j=N_k}\left(p_n^q\right)c_{m-1/2, i = 1/2}^{\a,k} + F_m^{\alpha,q,s} \right) \!\cdot\! g(p_n^q) \!\cdot\! R_\a T a_\alpha}
\end{align}
Rearranging equation \eqref{eq:nodal_pressure_nonideal_gas_appendix} leads to an explicit expression for the coefficient term $g(p_n^q)$:
\begin{align}\label{eq:g_coeff_as_fun_of_pres_nonideal_appendix}
	g(p_n^q) = \frac{p_n^q}{\sum_\a \left(\sum_{k\in\partial_+q} S_k \phi^{k}_{m,j=N_k}\left(p_n^q\right)c_{m-1/2, i = 1/2}^{\a,k} + F_m^{\alpha,q,s} \right) R_\a T (1 + a_\a p_n^q)}.
\end{align}
Substituting equation \eqref{eq:g_coeff_as_fun_of_pres_nonideal_appendix} into equation \eqref{eq:explicit_nodal_density_appendix_d} leads to an expression for nodal component densities after mixing at node $q$ with linear dependence on nodal pressure $p_n^q$ and explicit dependence on boundary flows and densities.  We write this dependence as $d_n^{\a,q} = d_n^{\a,q}(p_n^q)$, with form
\begin{align}\label{eq:explicit_nodal_density_nonideal_appendix}
	d_n^{\a,q}(p_n^q) = p_n^q \cdot \frac{\sum_{k\in\partial_+r} S_k \phi^{k}_{m,j=N_k}\left(p_n^q\right)c_{m-1/2, i = 1/2}^{\a,k} + F_m^{\alpha,q,s}}{\sum_\a \left(\sum_{k\in\partial_+r} S_k \phi^{k}_{m,j=N_k}\left(p_n^q\right)c_{m-1/2, i = 1/2}^{\a,k} + F_m^{\alpha,q,s} \right) R_\a T(1 + a_\a p_n^q)}.
\end{align}
We can then formulate the time update of boundary flows for each pipe $k\in\mathcal{E}$ based on the dependence on nodal pressure  $p^{r,n}$ as	
\begin{equation}\label{eq:boundary_flow_as_f_of_p_nonideal_appendix}
	\phi_{m,j=0}^{k}(p_n^q)
	\!=\! \left(\phi_{m-1,j=0}^{k} 	\!-\!
	\frac{\dt \lambda_k}{2D_k}
	\frac{\phi_{m-1,j=0}^{k}|\phi_{m-1,j=0}^{k}|}{\sum_\a d_{m - 1/2,i = 1/2}^{\a,k}} \!-\! \frac{\dt}{\dx}\cdot \frac{\sum_\a d_{m - 1/2,i = 1/2}^{\a,k} R_\a T}{1 - \sum_\a d_{m - 1/2,i = 1/2}^{\a,k} R_\a T a_\a} 
  \right) \!+\! \frac{\dt}{\dx} p_n^q.
\end{equation}
The time update for outgoing boundary densities for each gas component $\alpha$ for every pipe $k\in\mathcal{E}$ is obtained as
\begin{align}\label{eq:conservation_of_mass_discrete_nonideal_appendix}
	d_{m + 1/2,i = 1/2}^{\a,k} = d_{m - 1/2,i = 1/2}^{\a,k}  
	- 
	\frac{\dt}{\dx}
	\left( \phi^{k}_{m,j=1} \frac{ d^{\a,k}_{m - 1/2,i = 1/2}}{\sum_\beta d^{\beta,k}_{m - 1/2,i = 1/2}} - \phi^{k}_{m,j=0}(p_n^q) \frac{d_n^{\a,q}(p_n^q)}{\sum_\beta d_n^{\beta,q}(p_n^q)} \right).
\end{align}
As in the ideal gas modeling case, both updating expressions \eqref{eq:boundary_flow_as_f_of_p_nonideal_appendix} and \eqref{eq:conservation_of_mass_discrete_nonideal_appendix} have explicit dependence on data for previous time points on the staggered grid that are already known.  Equation \eqref{eq:conservation_of_mass_discrete_nonideal_appendix} is equivalent to equation \eqref{eq:conservation_of_mass_discrete_ideal_appendix}, but with equation \eqref{eq:explicit_nodal_density_nonideal_appendix} used for computing component nodal densities after mixing.

\subsection{Input Flow Monitoring}  \label{sec:append_monitoring}

Here we provide details on derivation of a local nodal monitoring policy that will guarantee that upper bounds on hydrogen mass fraction are satisfied in real time.  In the motivating scenario, there is a gas mixture flowing at a constant rate into a given node $q$ with known individual gas concentrations.  We suppose that the composition of this injection flow is different from the composition of gases flowing into the node from incoming network pipes $k\in\partial^+q$, and the nodal injection thus alters the composition of the flow leaving the node to outgoing pipes $k\in\partial^-q$.  If at some point, the mass fraction $c_m^{\alpha,q}$ at node $q$ of a constituent gas $\alpha$ approaches a safety limit $c_{\max}^{\alpha,q}$, the monitoring and control policy goal is to compute a gas injection rate at node $q$ that ensures that this mass fraction limit is not exceeded.

Recall equation \eqref{eq:conservation_of_momentum_discrete_flux_pipe_ideal_appendix} (or \eqref{eq:conservation_of_momentum_discrete_flux_pipe_nonideal_appendix}) for pipe boundary flow updating as a function of nodal pressure and data at the previous time point, 
\begin{equation} \label{eq:conservation_of_momentum_discrete_flux_pipe_appendix}
    \phi^{k}_{m,j=0}(p^q_n)
	= \! \left(\!\phi_{m-1,j=0}^{k}
	-
	\frac{\dt \lambda_k}{2D_k}
	\frac{\phi_{m-1,j=0}^{k}|\phi_{m-1,j=0}^{k}|}{\sum_\a d^{\a,k}_{m - 1/2,i = 1/2}} \!-\! \frac{\dt}{\dx} p_{n,i=1/2}^{k} \!\right) \!+\! \frac{\dt}{\dx} \mu_n^k p^q_n.
\end{equation}
and equation \eqref{eq:flow_balance_scalar_coef1_appendix} for computing pressure at node $q$,
\begin{equation} \label{eq:flow_balance_scalar_coef1_injection_only_appendix}
	 p_n^{q} \!=\! \frac{\dx}{\dt}\!\cdot\! \frac{ F^{q,s}_n - \sum_{k\in\partial q} S_k  \left(\phi_{m-1,j=0}^{k} - \frac{\dt \lambda_k}{2D_k} \cdot \frac{\phi_{m-1,j=0}^{k}|\phi_{m-1,j=0}^{k}|}{\sum_\a d_{m - 1/2,i = 1/2}^{\a,k}} - \frac{\dt}{\dx} p_{n,i=1/2}^{k} \right)}{\sum_{k\in\partial q} S_k\mu_n^k},
\end{equation}
where $\mu_n^k\equiv 1$ for all pipes $k$ with no compressors.  For the purpose of this derivation, we suppose that there is only injection and no withdrawal at node $q$ so that $F^{q,d}_n\equiv 0$ in equation \eqref{eq:flow_balance_scalar_coef1_injection_only_appendix}.
Observe in equation \eqref{eq:flow_balance_scalar_coef1_injection_only_appendix} that because $p^q_n$ depends on the total nodal gas injection flow $F_n^{q,s}$, then $p^q_n$ can be written as $p^q_n(F_n^{q,s})$ and we can in turn rewrite $\phi^{k}_{m,j=0}(p^q_n)$ as $\phi^{k}_{m,j=0}(F_n^{q,s})$:
\begin{align} 
    \phi^{k}_{m,j=0}(F_n^{q,s})
	&= \! \left(\!\phi_{m-1,j=0}^{k}
	-
	\frac{\dt \lambda_k}{2D_k}
	\frac{\phi_{m-1,j=0}^{k}|\phi_{m-1,j=0}^{k}|}{\sum_\a d^{\a,k}_{m - 1/2,i = 1/2}} \!-\! \frac{\dt}{\dx} p_{n,i=1/2}^{k} \!\right) \nonumber \\
    & \, +\! \mu_n^k \cdot\! \frac{ F^{q,s}_n- \sum_{k\in\partial q} S_k  \left(\phi_{m-1,j=0}^{k} - \frac{\dt \lambda_k}{2D_k} \cdot \frac{\phi_{m-1,j=0}^{k}|\phi_{m-1,j=0}^{k}|}{\sum_\a d_{m - 1/2,i = 1/2}^{\a,k}} - \frac{\dt}{\dx} p_{n,i=1/2}^{k} \right)}{\sum_{r\in\partial q} S_r\mu_n^r}.  \label{eq:boundary_flow_depending_on_injection_appendix}
\end{align}
To simplify subsequent derivation, we re-write equation \eqref{eq:boundary_flow_depending_on_injection_appendix} using a shorthand notation as
\begin{equation} \label{eq:boundary_flow_depending_on_injection_shorthand_appendix}
         \phi^{k}_{m,j=0}(F_n^{q,s}) = \frac{\mu_n^k F_n^{q,s}}{\sum_r \mu_n^r S_r} + \Theta_m^{q,k} + \mu_n^k \cdot \Upsilon_m^q.
\end{equation}
where
\begin{align*}
    & \Theta_m^{q,k} = \left(\!\phi_{m-1,j=0}^{k}
	-
	\frac{\dt \lambda_k}{2D_k}
	\frac{\phi_{m-1,j=0}^{k}|\phi_{m-1,j=0}^{k}|}{\sum_\a d^{\a,k}_{m - 1/2,i = 1/2}} \!-\! \frac{\dt}{\dx} p_{n,i=1/2}^{k} \!\right), \quad \forall k \in \partial q, \\
	& \Upsilon_m^q = \frac{ \sum_{k\in\partial q} S_k  \left( \frac{\dt \lambda_k}{2D_k} \cdot \frac{\phi_{m-1,j=0}^{k}|\phi_{m-1,j=0}^{k}|}{\sum_\a d_{m - 1/2,i = 1/2}^{\a,k}} - \frac{\dt}{\dx} p_{n,i=1/2}^{k} - \phi_{m-1,j=0}^{k}\right)}{\sum_{k\in\partial q} S_k\mu_n^k}.
\end{align*}
Denoting the maximum mass fraction limit for gas $\alpha$ at node $q$ as $c^{\alpha,q}_{\max}$, supposing that we wish to adjust mass injection $F_n^{q,s}$ to maintain the mass fraction at that limit, we can apply equation \eqref{eq:explicit_nodal_density2_appendix_c} to examine the required condition:
\begin{align} \label{eq:monitoring_mass_fraction_setting_appendix}
    c^{\alpha,q}_{\max}= c_m^{\alpha,q} & = \frac{\sum_{k\in\partial_+r} S_k \phi^{k}_{m,j=N_k}(F_n^{q,s})c_{m-1/2, i = 1/2}^{\a,k} + F_m^{q,s}\cdot c_m^{\alpha,s} }{\sum_{k\in\partial_-r}  S_k \phi^{k}_{m,j=0}(F_n^{q,s}) },  
\end{align}
where we have written boundary flows $\phi^{k}_{m,j=0}$ as functions of the total injection flow $F_n^{q,s}$, and decomposed the component inflow $F_n^{\alpha,q,s}$ as total inflow $F_m^{q,s}$ times component mass fraction $c_m^{\alpha,s}$. Rearranging and substituting for boundary flows using equation \eqref{eq:boundary_flow_depending_on_injection_shorthand_appendix} leads to
\begin{align}
& \sum_{k\in\partial_-r}  S_k \left( \frac{\mu_n^k F_n^{q,s}}{\sum_r \mu_n^r S_r} + \Theta_m^{q,k} + \mu_n^k \cdot \Upsilon_m^q \right) c^{\alpha,q}_{\max} = \nonumber \\
& \qquad \sum_{k\in\partial_+r} S_k \left( \frac{\mu_n^k F_n^{q,s}}{\sum_r \mu_n^r S_r} + \Theta_m^{q,k} + \mu_n^k \cdot \Upsilon_m^q \right) c_{m-1/2, i = 1/2}^{\a,k} + F_m^{q,s}\cdot c_m^{\alpha,s}.
\end{align}
Collecting terms and solving for $F_n^{q,s}$ leads to a formula for the total injection at node $q$ that maintains the nodal mass fraction after mixing at node $q$ at the maximum allowable level $c^{\alpha,q}_{\max}$. We write this expression as a function $F_n^{q,s}(c^{\alpha,q}_{\max})$:
\begin{equation} \label{eq:max_injection_at_mass_fraction_limit_appendix}
F_n^{q,s}(c^{\alpha,q}_{\max}) \!=\! \frac{c_{m-1/2, i = 1/2}^{\a,k} \sum_{k\in\partial_+r} S_k \left(\Theta_m^{q,k} + \mu_n^k \!\cdot\! \Upsilon_m^q \right) \!-\! c^{\alpha,q}_{\max} \sum_{k\in\partial_-r}  S_k \left( \Theta_m^{q,k} + \mu_n^k \!\cdot\! \Upsilon_m^q \right) }{  c^{\alpha,q}_{\max} \sum_{k\in\partial_-r}  S_k \left( \frac{\mu_n^k}{\sum_r \mu_n^r S_r} \right) - c_{m-1/2, i = 1/2}^{\a,k} \sum_{k\in\partial_+r} S_k \left( \frac{\mu_n^k}{\sum_r \mu_n^r S_r} \right) - c_m^{\alpha,s}}.
\end{equation}
Specifically, we may suppose that $c_m^{\alpha,s}\equiv 1$ where $\alpha$ denotes hydrogen mass fraction, and $c^{\alpha,q}_{\max}$ is the maximum allowable hydrogen mass fraction inside the pipeline system.

\subsection{Compressibility factors} \label{appendix:compressibility}

We summarize the approaches that we use to obtain a linear approximations to the compressibility factors $Z_{H2}$ and $Z_{NG}$ of hydrogen and natural gas using tabulated experimental data from previous studies.  The data used to develop the approximations for hydrogen and natural gas and  are available in the sources \cite{mihara1977comph2} and \cite{kareem2016compng}, respectively. We propose the use of these linear approximations because of the nearly linear dependence of both compressibility factors on pressure, as seen in the experimental measurements.  For simplicity, as also done in the numerical approach developed above, we suppose that temperature is constant at $\mathbf{T} = 298.15^{\circ}\mbox{K}$.

\subsubsection{Hydrogen Compressibility}

Experimental measurements from a previous study \cite{mihara1977comph2} can be compiled to tabulate the dependence of the compressibility factor of hydrogen on pressure (atm), as shown in table \ref{tab:hydrogen_compressibility}.

\begin{table}[ht!]
    \centering
    \begin{tabular}{|| c | c | c || c | c | c ||} 
    \hline
    pressure   & compressibility & $a_h$ & pressure  & compressibility & $a_h$\\
     (atm) &  $Z_{H2}$ &  &  (atm) &  $Z_{H2}$ & \\
    [0.5ex] 
    \hline\hline
    83.731 & 1.0503 &  5.999 $\times 10^{-8}$ & 13.478 & 1.0078 &   5.780 $\times 10^{-8}$ \\ 
    \hline
    69.748 & 1.0417 &   5.971 $\times 10^{-8}$ & 11.306 & 1.0066 &  5.830 $\times 10^{-8}$ \\
    \hline
    52.613 & 1.0312 &  5.922 $\times 10^{-8}$ & 8.6021 & 1.0050 &  5.805 $\times 10^{-8}$ \\
    \hline
    43.964 & 1.0260 &  5.906 $\times 10^{-8}$ & 7.2187 & 1.0042 &  5.811 $\times 10^{-8}$ \\
    \hline
    33.291 & 1.0196 &  5.880 $\times 10^{-8}$ & 5.4954 & 1.0032 &  5.815 $\times 10^{-8}$ \\
    \hline
    27.872 & 1.0163 &  5.840 $\times 10^{-8}$ & 4.6129 & 1.0027 &  5.845 $\times 10^{-8}$ \\ 
    \hline
    21.155 & 1.0123 &  5.806 $\times 10^{-8}$ & 3.5129 & 1.0021 &  5.970 $\times 10^{-8}$ \\
    \hline
    17.732 & 1.0103 & 5.801 $\times 10^{-8}$ \\
    \hline
    \end{tabular}
    \caption{Experimental compressibility factors for hydrogen \cite{mihara1977comph2} at fixed temperature $\mathbf{T} = 298.15 ^{\circ}\mbox{K}$.}
    \label{tab:hydrogen_compressibility}
\end{table}

Supposing that the linear approximation of $Z_{H2}$ as a function of pressure at constant temperature $\mathbf{T} = 298.15 ^{\circ}\mbox{K}$ is given as $Z_{H2}(p) = 1 + a_{h}(p) \cdot p$, then $a_{h}$ values can be computed, e.g., as $5.97\times10^{-8}=(1.0021-1)/(3.5129 \cdot 101,325)$, from the pressure and compressibility values listed in table \ref{tab:hydrogen_compressibility}, where $101,325$ is the conversion factor between Pa and atm.  Taking the mean of the $a_h$ results from table \ref{tab:hydrogen_compressibility}, we estimate the linear coefficient  $a_{h}$  in the approximation of the compressibility factor for hydrogen as follows:
\begin{align*}
    & Z(p) = 1 + a_{h}(p) \cdot p, 
    \qquad 
    a_{h} = \frac{Z(p) - 1}{p}
    %\\  & a_{h} 
    \approx 
     %\frac{1.0021 - 1}{3.5129 \cdot 101,325} \approx 
    5.865 \times 10^{-8}.
    %0.0000000059,
\end{align*}

\subsubsection{Natural Gas Compressibility}

An appropriate linear approximation to the compressibility factor for natural gas can be found by using the results in \cite{kareem2016compng}. Pseudo-reduced temperature and pressure are defined in \cite{dranchuk1971} as the ratio of temperature (psi) and pressure ($^{\circ}\mbox{R}$) to the pseudo-critical temperature ($^{\circ}\mbox{R}$) and pressure (psi) of natural gas, respectively:
\begin{align}
    T_{pr} = \frac{T}{T_{pc}}, \quad P_{pr} = \frac{P}{P_{pc}}.
\end{align}
Pseudo-critical temperature and pressure can be written in terms of specific gravity \cite{sutton1985compressibility}:
\begin{align}\label{eq:Pseudocritical_T_and_P}
    T_{pc} &= 169.2 + 349.5 \gamma_g - 74.0 \gamma_g^2,\\
    P_{pc} &= 756.8 - 131.07 \gamma_g - 3.6 \gamma_g^2.
\end{align}
Fixing (for simplicity) the specific gravity of natural gas at $\gamma_g = 0.7$, we apply equation \eqref{eq:Pseudocritical_T_and_P} to get pseudo-critical temperature and pressures for natural gas: 
\begin{align*}
    T_{pc}^{ng} &= 169.2 + 349.5 \cdot 0.7 - 74.0 \cdot 0.49 = 169.2 + 244.65 - 36.26 = 377.59,\\
    P_{pc}^{ng} &= 756.8 - 131.07\cdot 0.7 - 3.6 \cdot 0.49 = 756.8 - 91.749 - 1.764 = 663.287.
\end{align*}
We then obtain the pseudo-reduced temperature using fixed temperature $T = 298.15^{\circ}\mbox{K} = 536.67 ^{\circ}\mbox{R}$ and \eqref{eq:Pseudocritical_T_and_P}:
\begin{align*}
    T_{pr}^{ng} = \frac{T}{T_{pc}^{ng}} = \frac{536.67}{377.59} = 1.4213.
\end{align*}
The form of the linear approximation of $Z_{NG}$ that we seek can be written in terms of pseudo-reduced pressure and then adapted to appear in terms of pressure:
\begin{align}\label{eq:compressibility_pressure_pseudo-reduced}
    Z(P_{pr}) &= 1 + a_1 P_{pr},\\
    Z(p)      &= 1 + a_{ng} p = 1 + \frac{a_1}{P_{pc}^{ng}} p.
\end{align}
Using data from Fig. 1 in \cite{kareem2016compng} that shows a plot of experimental measurements of the $Z_{NG}$, we have $Z \approx 0.7666$, $P_{pr} = 2$. Using equation \eqref{eq:compressibility_pressure_pseudo-reduced} yields
\begin{align*}
    a_1    &= \frac{0.7666 - 1}{2} \approx -0.1167,\\
    a_{ng} &= \frac{-0.1167}{663.287} \cdot 0.000145 \approx -0.25 \times 10^{-7},
    %0.000000025,
\end{align*}
where $0.000145$ is the conversion factor from Pa to psi.

\section*{Acknowledgements}  The authors are grateful to Falk M. Hante and Saif R. Kazi for numerous helpful discussions, and to Luke S. Baker for sharing the lumped-element model simulation used as a comparison benchmark. This study was supported by the U.S. Department of Energy's Advanced Grid Modeling (AGM) project ``Dynamical Modeling, Estimation, and Optimal Control of Electrical Grid-Natural Gas Transmission Systems'', as well as LANL Laboratory Directed R\&D project ``Efficient Multi-scale Modeling of Clean Hydrogen Blending in Large Natural Gas Pipelines to Reduce Carbon Emissions''. Research conducted at Los Alamos National Laboratory is done under the auspices of the National Nuclear Security Administration of the U.S. Department of Energy under Contract No. 89233218CNA000001.  Report No. LA-UR-24-22948.

\bibliography{main.bib}

\end{document}